\documentclass[11pt]{article}
\usepackage{fullpage}
\usepackage{graphicx}
\usepackage{amsmath}
\usepackage{amsfonts}
\usepackage{amssymb}
\usepackage{epstopdf}
\usepackage{amsthm}
\usepackage{url}
\usepackage{cite}
\usepackage{subfigure}

\graphicspath{{./figures/}}

\begin{document}
\title{A Modulus-Squared Dirichlet Boundary Condition for
Time-Dependent Complex Partial Differential Equations and
its Application to the Nonlinear Schr{\"o}dinger Equation}
\author{
R. M. Caplan\footnote{Corresponding author.
Present address: Predictive Science Inc.
9990 Mesa Rim Rd, Suite 170, San Diego, CA 92121.
 email: caplanr@predsci.com, phone: 858-225-2314,
\texttt{URL}: http://www.predsci.com}\,
and
R. Carretero-Gonz{\'a}lez\\[1.0ex]
Nonlinear Dynamical System Group\footnote{\texttt{URL}: http://nlds.sdsu.edu},
Computational Science Research Center\footnote{\texttt{URL}: http://www.csrc.sdsu.edu},
and \\
Department of Mathematics and Statistics,
San Diego State University,\\
San Diego, California 92182-7720, USA\\
}

\date{}
\maketitle

\begin{abstract}
An easy to implement modulus-squared Dirichlet (MSD) boundary condition is formulated for numerical simulations of time-dependent complex partial differential equations in multidimensional settings.  The MSD boundary condition approximates a constant modulus-square value of the solution at the boundaries and is defined as
\[
\left.\dfrac{\partial \Psi}{\partial t}\right|_b \approx i\,\mbox{Im}\left[  \dfrac{1}{\Psi_{b-1}}  \left.\dfrac{\partial \Psi}{\partial t}\right|_{b-1}  \right]\,\Psi_b,
\]
where $\Psi$ is the complex field and
the subscripts $b$ and $b-1$ refer to a boundary point and the closest interior point to the boundary respectively.
Application of the MSD boundary condition to simulations of the nonlinear Schr{\"o}dinger equation is shown, and numerical simulations are performed to demonstrate its usefulness and advantages over other simple boundary conditions.
\end{abstract}

\section{Introduction}
\label{s:intro}
When utilizing numerical methods to approximate the solutions to time-dependent partial differential equations (PDEs), proper handling of boundary conditions can be quite challenging. Sometimes, an otherwise stable numerical scheme will become unstable depending on how the boundary conditions are computed \cite{FD_PDE_BOOK}. In addition, high-order schemes can degrade in accuracy to lower-order when using boundary conditions which are not compatible with the high-order accuracy \cite{FD_PS_BOOK}.  Proper handling of boundary conditions in higher-order schemes, especially in high-order compact schemes can be even more of a challenge \cite{MSD_BC_HOC_hyperbolic,MSD_BC_HOC}.

Often, researchers will forgo a complicated boundary condition implementation and instead use tried-and-true boundary condition techniques which are very simple yet provide acceptable results. One of the most common is the use of Dirichlet boundary conditions when simulating solutions which decay towards zero at infinity, and where most of the dynamics (or `action') is expected to remain in the central regions of the computational grid.  Another simple method in such cases is to use periodic boundary conditions.

Infinite-domain problems involving PDEs whose function values are complex cannot, in general, make use of numerical Dirichlet or periodic boundary conditions because of the oscillation of the real and imaginary parts of the function due to the intrinsic frequency of the system.  (In cases when both the real and imaginary parts of the solution converge to a constant at infinity, such boundary conditions \emph{can} be used).  In many solutions (for examples, see Sec.~\ref{s:num}), the real and imaginary parts of the solution do not converge, but rather, the modulus-squared of the solution converges to a constant value at infinity.  In such a case, Dirichlet boundary conditions fail and it may not be possible to use periodic boundary conditions since the solution will instantaneously `feel' its periodic counterpart through its phase which can cause dynamical effects which would not happen in a true infinite background.  Additionally, for some solutions, a large phase jump can exist between one boundary of the domain and the opposite boundary in which case using periodic boundary conditions artificially connects the phase jump resulting in an incorrect solution.  It is therefore desirable in such cases to find a simple alternative boundary condition.

In this paper, we present a simple way to simulate a modulus-squared Dirichlet (MSD) boundary condition in time-dependent complex-valued PDEs which keeps the modulus-squared of the solution at the boundaries constant.  The MSD boundary condition is very easy to implement and eliminates the need for overly large grids or expensive and complicated boundary conditions for many problems.

A very common time-dependent complex-valued PDE used in a wide range of applications is the nonlinear {Schr{\"o}dinger} equation (NLSE).  The NLSE is a universal model describing the evolution and propagation of complex field envelopes in nonlinear dispersive media.  As such, it is used to describe many physical systems including the evolution of water waves, nonlinear optics, thermodynamic pulses, nonlinear waves in fluid dynamics, waves in semiconductors, and the mean-field dynamics of Bose-Einstein condensates (BECs) \cite{NLSE_nlpdebook, NLSE_selffocus,NLSE_OPTVA,BEC_RCbook}.  In systems such as optics and BECs, the modulus-squared of the solution (referred to as the `wavefunction') represents the observable (intensity of light and atomic density respectively).  In this situation, often the dynamics of `dark' structures (dark solitons, vortices, vortex lines, vortex rings, etc.) which reside inside the medium are studied.  They are coherent structures of very low (or zero) central density which exist inside the bulk of the medium.  The most basic form of the structures can be examined by assuming an infinite-extent domain, in which case the solutions exist within an infinite constant-density background.  Such a situation is very well-suited for the use of a MSD boundary condition.  As such, in Sec.~\ref{s:num}, we use simulations of dark coherent structures in the NLSE to test the MSD boundary condition.

In the specific case of the linear and nonlinear {Schr{\"o}dinger} equations,
many sophisticated boundary conditions have been developed which simulate
transparent or artificial boundaries.
%
For an extensive review on this topic
we refer the interested reader to the Reviews \cite{MSD_NLSE_TRANS_REV,Bao_CPC_13}.
For the linear {Schr{\"o}dinger} equation, such boundary conditions include continuous, pole condition, temporally discrete, spatially discrete, fully discrete, as well as others in one dimension, while in two dimensions, continuous transparent boundary conditions have been formulated which make use of the Sommerfeld radiation condition \cite{MSD_RADBCmultidim}.
Work has also been done in formulating similar conditions for the one-dimensional NLSE (see Refs.~\cite{MSD_NLSE_TRANS_REV,Bao_CPC_13,Bao_KRM_13} and references therein for a complete review).  Most of the aforementioned boundary conditions focus on eliminating reflections off the boundaries when studying dynamics of solutions which tend to zero at infinity.
For example, the method of perfectly matched layers
\cite{Berenger_JCP_94,Neuhauser_BIT_90,Zheng_JCP_07} has been used
to artificially damp traveling waves in the domain without
modifying their dynamics.
Another popular method to deal with boundary conditions is to
implement absorbing boundary conditions
\cite{Antoine_SJNA_06,Antoine_JCP_09,Antoine_LP_11,Antoine_SJSC_11,Antoine_MMMAS_12,Antoine_NM_13}.
While these boundary conditions would be valid for the study of bright structures (a localized solution on a zero background) in the NLSE, they are not necessarily appropriate for dark structures (hole or ``bubble'') which exist in a constant-density, non-zero, background.
Furthermore, transparent or artificial boundary conditions can be very complicated to implement.  A relevant exception to this is the Sommerfeld radiation condition described in Ref.~\cite{MSD_RADBCmultidim}, which is relatively easy to implement and has been previously successful in simulating a constant-density background in the NLSE \cite{VR_aroundball_usesradBC_mbc7}.  However, since it requires additional computations at the boundaries, and there has been some discussion debating its usefulness in some cases (see Ref.~\cite{MSD_RADBAD}), the MSD boundary condition formulated in this paper remains a useful alternative approach to approximate a constant-density background in the NLSE.

The paper is organized as follows.  In Sec.~\ref{s:formation} we formulate the MSD boundary condition for general time-dependent complex-valued PDEs.  Then, in Sec.~\ref{s:nlseapp}, we apply the MSD boundary condition to the NLSE.  In Sec.~\ref{s:num} we numerically demonstrate the usefulness of the MSD boundary condition for simulating the one-, two-, and three-dimensional NLSE using a Runge-Kutta finite-difference scheme and comment on the stability of the MSD boundary condition. Finally, we conclude in Sec.~\ref{s:con}.

\section{Formulation of the MSD Boundary Condition}
\label{s:formation}
\subsection{Notation}
\label{s:note}
To formulate a MSD boundary condition we first introduce the notation that $\Psi_b$ describes the grid-points on a boundary and $\Psi_{b-1}$ describes the grid-points one cell in from the boundary in the normal direction to the boundary points.  The real part of $\Psi$ is denoted as $\Psi^R$ and the imaginary part as $\Psi^I$.  The short-hand notation of $\Psi_{\alpha}$ is used to denote the first derivative with respect to $\alpha$ ($\partial \Psi / \partial \alpha$).
\subsection{Derivation of the MSD boundary Condition}
\label{s:derive}
We start by stating the condition, that for all times, the modulus-squared of the function $\Psi$ at the boundaries ($\Psi_b$) is equal to a constant, positive, real value $B$:
\begin{equation}
\label{msdassumed}
|\Psi_b|^2 = B.
\end{equation}
This is equivalent to stating that the solution at the boundary has the form
\begin{equation}
\label{steadystatepsi}
\Psi_b = \sqrt{B} \, \mbox{exp}\left[-i\;\hat \Omega_b(\vec r,t)\right],
\end{equation}
where $\hat \Omega(\vec r,t)\in \Re$ is the phase which can be a function of space and time.  Taking the natural logarithm of Eq.~(\ref{steadystatepsi}) yields
\[
\mbox{ln}(\Psi_b) = \mbox{ln}\left(\sqrt{B}\right) - i\,\hat \Omega_b(\vec r,t),
\]
which, after taking a temporal derivative, yields
\begin{equation}
\label{imagness}
\frac{\partial}{\partial t}\,\mbox{ln}(\Psi_b) = - i\,\frac{\partial}{\partial t}\,\hat \Omega_b(\vec r,t).
\end{equation}
We now apply a spatial derivative in the normal direction to the boundary (denoted here by $x$) on each side yielding the continuous representation of the MSD boundary condition:
\begin{equation}
\label{msd_cont}
\frac{\partial}{\partial x}\,\left(\frac{\Psi_{t,b}}{\Psi_b}\right) = -i\,\hat \Omega_{xt,b}(\vec r,t).
\end{equation}
Discretizing the left-hand side of Eq.~(\ref{msd_cont}) using first-order one-sided differencing and multiplying both sides by $h_x$ (the grid spacing in the normal direction to the boundary) yields
\begin{equation}
\label{msddisc1}
\frac{\Psi_{t,b}}{\Psi_b} - \frac{\Psi_{t,b-1}}{\Psi_{b-1}}
 = -i\,h_x\,\hat \Omega_{xt,b}(\vec r,t) + O(h_x^2).
\end{equation}
We note that if the phase of the solution at the boundary is
a linear combination of space and time dependent functions, namely
$\hat \Omega(\vec r,t) = \alpha_1 f(\vec r) + \alpha_2 g(t)$ for any real
$\alpha_1$ and $\alpha_2$, then the $\hat \Omega_{xt,b}$ term of Eq.~(\ref{msddisc1}) 
becomes exactly zero, while in the general case (where there could be a combined spatio-temporal 
phase at the boundary that yields a non-zero $\hat \Omega(\vec r,t)$ term), 
the right hand side of Eq.~(\ref{msddisc1}) remains $O(h_x)$.

We approximate the right hand side of Eq.~(\ref{msddisc1}) as zero, noting that
for practical applications such as the simulation of coherent structures which are
sufficiently far away from the boundaries (see Sec.~\ref{s:nlseapp}), the phase
is expected to be a {\em slowly} varying function ---at most a linear function---
of time and space at the boundaries.  
Therefore,
after rearranging Eq.~(\ref{msddisc1}), we arrive at the discrete form of the MSD boundary condition:
\begin{equation}
\label{msd}
\Psi_{t,b} \approx \frac{\Psi_{t,b-1}}{\Psi_{b-1}}\,\Psi_b,
\end{equation}
where the value of $\Psi_{t,b-1}$ is computed using the interior numerical scheme being implemented.  

It is interesting to note that the MSD boundary condition of Eq.~(\ref{msd}) is very similar in form to the Sommerfeld radiation boundary condition mentioned in Ref.~\cite{MSD_RADBCmultidim} where the time-derivative at the boundary is computed as
\[
\Psi_{t,b} \approx -\frac{\Psi_{t,b-1}}{\Psi_{x,b-1}}\,\Psi_{x,b},
\]
where $x$ is the direction normal to the boundary.




It is important to note that from Eq.~(\ref{imagness}), $\Psi_{t,b}/\Psi_b$ should be purely imaginary, and therefore from Eq.~(\ref{msd}), $\Psi_{t,b-1}/\Psi_{b-1}$ should be purely imaginary as well.  However, due to numerical errors, computing $\Psi_{t,b-1}/\Psi_{b-1}$ may introduce a small real part which would cause the solution at the boundaries to undergo eventual exponential growth (which we have observed in simulations not reported here).  In order to ensure that $\Psi_{t,b-1}/\Psi_{b-1}$ is purely imaginary, we modify the MSD boundary condition of Eq.~(\ref{msd}) to be
\begin{equation}
\label{msdfixed}
\Psi_{t,b} \approx i\,\mbox{Im}\left[\frac{\Psi_{t,b-1}}{\Psi_{b-1}}\right]\,\Psi_b.
\end{equation}

When using the MSD boundary condition in programming environments that do not intrinsically handle complex variables, Eq.~(\ref{msdfixed}) must be explicitly split into its real and imaginary parts.  We begin by expanding the unmodified MSD boundary condition of Eq.~(\ref{msd}) into its real and imaginary parts yielding
\[
\Psi_{t,b}^R + i\,\Psi_{t,b}^I \approx \frac{\Psi_{t,b-1}^R + i\,\Psi_{t,b-1}^I}{\Psi_{b-1}^R + i\,\Psi_{b-1}^I}\left(\Psi_{b}^R + i\,\Psi_{b}^I\right).
\]
which leads to
\begin{alignat}{2}
\label{msdri2}
\Psi_{t,b}^R + i\,\Psi_{t,b}^I \approx &\left[\frac{\Psi_{t,b-1}^R\,\Psi_{b-1}^R + \Psi_{t,b-1}^I\,\Psi_{b-1}^I}{\left(\Psi^R_{b-1}\right)^2 + \left(\Psi^I_{b-1}\right)^2} \right. \\
&\left. + i\,\left(\frac{\Psi_{t,b-1}^I\,\Psi_{b-1}^R - \Psi_{t,b-1}^R\,\Psi_{b-1}^I}{\left(\Psi^R_{b-1}\right)^2 + \left(\Psi^I_{b-1}\right)^2}\right) \right]\,\left(\Psi_{b}^R + i\,\Psi_{b}^I\right). \notag
\end{alignat}
We note that with the MSD assumption $|\Psi_b|^2 = B$, we have that
 \begin{equation}
 \label{msdode1}
 \frac{\partial}{\partial t}\,|\Psi_b|^2 = \frac{\partial}{\partial t}\,\left(\Psi_b^R\right)^2 + \frac{\partial}{\partial t}\,\left(\Psi_b^I\right)^2 = 0,
 \end{equation}
 which is equivalent to
 \begin{equation}
 \label{msdode2}
 \Psi_b^R\,\frac{\partial\Psi_b^R}{\partial t} =-\Psi_b^I\,\frac{\partial\Psi_b^I}{\partial t}.
 \end{equation}
Therefore, we see that the first term in the brackets of Eq.~(\ref{msdri2}) would be equal to zero if evaluated at the boundary.  If it is assumed that the interior point is similar, the term can be dropped.  Dropping the term is equivalent to the numerical fix used in Eq.~(\ref{msdfixed}).  Simplifying Eq.~(\ref{msdri2}) with this in mind yields the separated MSD boundary condition
\begin{alignat}{4}
\label{ab}
\Psi^R_{t,b} &= -\tilde \Omega \, \Psi^I_{b},\\
\Psi^I_{t,b} &=  \tilde \Omega \, \Psi^R_{b}, \notag
\end{alignat}
where
\begin{equation}
\label{msdom}
\tilde \Omega = \frac{\Psi^I_{t,b-1}\,\Psi^R_{b-1} - \Psi^R_{t,b-1}\,\Psi^I_{b-1}}{\left(\Psi^R_{b-1}\right)^2 + \left(\Psi^I_{b-1}\right)^2}.
\end{equation}

The MSD boundary condition of Eqs.~(\ref{msdfixed}) and (\ref{ab}) is given for the temporal derivative of the boundary point.  This is ideally suited for Runge-Kutta type solvers, as the right hand side of the PDE ($\Psi_t$) is evaluated and used for the time-stepping \cite{RK4}.  For other methodologies, or in situations where the boundary value of the spatial derivatives is needed, the MSD boundary condition can be inserted into the governing equation to formulate the required boundary conditions.  An example of this is shown in Sec.~\ref{s:nlseapp} for the NLSE.

In situations where the sequential computation of the internal scheme followed by the boundary condition computations is not appropriate (for example, implicit finite-difference schemes) one can substitute the internal scheme of $\Psi_{t,b-1}$ into Eq.~(\ref{msdfixed}) and implement the boundary condition concurrently with the internal scheme (this would add extra computational steps as $\Psi_{t,b-1}$ would essentially be computed twice).

If it happens that $\Psi_{b-1} = 0$, the MSD boundary condition of Eq.~(\ref{msd}) has a singularity.  However, in most situations, $\Psi_{b-1}$ only takes a zero value when the solution is tending toward zero at the boundary, in which case, the standard Dirichlet boundary condition of $\Psi_b = 0$ can be used instead of the MSD boundary condition.  Alternatively, one can numerically check the values of $\Psi_{b-1}$, and use an alternative boundary condition in the case that $\Psi_b = 0$.










\subsection{Key Features of the MSD Boundary Condition}
\label{s:features}
The following are a few key features of the MSD boundary condition:
\medskip

\begin{itemize}
\item The MSD boundary condition does not depend on the specific PDE being simulated, only that it is time-dependent and complex-valued, nor does it depend on the numerical scheme used to compute the internal points of the grid.
\item The MSD boundary condition can be used in multidimensional settings without modification since each boundary point only uses one interior point in the normal direction to the boundary.
\item The MSD boundary condition does not depend on the size of the grid spacing, and therefore does not need to be altered when using non-uniform grid spacing (for example, in adaptive mesh refinement applications).
\item The MSD boundary condition can be considered compact, in that after computing the internal scheme, the boundary values only depend on their nearest neighboring grid point.  This allows the MSD boundary condition to be easily implemented in parallel environments.
\item In general, the MSD boundary condition is extremely easy to implement. This makes it an attractive alternative to more complicated boundary condition methodologies.
\end{itemize}

\section{Application of the MSD boundary condition to the NLSE}
\label{s:nlseapp}
Here we show an implementation of the MSD boundary condition for the NLSE.  The NLSE with an external potential can be given in general form as
\begin{equation}
\label{NLSE}
i\,\frac{\partial \Psi}{\partial t} + a\,\nabla^2\Psi - V({\bf r})\,\Psi +s\,|\Psi|^2\,\Psi = 0,
\end{equation}
where $\Psi$ is the wavefunction, $V({\bf r})$ is an external potential, and $a$ and $s$ are parameters determined by the physical system being modeled.

As mentioned in Sec.~\ref{s:derive}, since the MSD boundary condition of Eq.~(\ref{msdfixed}) is given in terms of the temporal derivative, it is well suited for Runge-Kutta schemes, in which case it can be applied directly.  However, other numerical schemes require boundary conditions on the Laplacian operator itself. In the case of the NLSE, this can be worked out by inserting the NLSE of Eq.~(\ref{NLSE}) into Eq.~(\ref{msdfixed}) yielding
\begin{equation}
\label{MSDlap1}
\nabla^2\Psi_b \approx \left[ \mbox{Im}\left(i\, \frac{\nabla^2\Psi_{b-1}}{\Psi_{b-1}}\right) + \frac{1}{a}\,\left(N_{b-1} - N_b\right)\right]\Psi_b,
\end{equation}
where
\begin{equation}
\label{nbnb1}
N_b = s\,|\Psi_b|^2 - V_b, \qquad N_{b-1} = s\,|\Psi_{b-1}|^2 - V_{b-1}.
\end{equation}
Splitting Eq.~(\ref{MSDlap1}) into real and imaginary parts yields
\begin{alignat}{3}
\nabla^2 \Psi^R_{b} &\approx \left[A + \frac{1}{a}\left(N_{b-1} - N_b\right)\right]\,\Psi^R_b, \\
\nabla^2 \Psi^I_{b} &\approx \left[A + \frac{1}{a}\left(N_{b-1} - N_b\right)\right]\,\Psi^I_b, \notag
\end{alignat}
where
\begin{equation}
\label{A}
A =  \frac{\nabla^2\Psi^R_{b-1}\,\Psi^R_{b-1} + \nabla^2\Psi^I_{b-1}\,\Psi^I_{b-1}}{\left(\Psi^R_{b-1}\right)^2 + \left(\Psi^I_{b-1}\right)^2},
\end{equation}
and $N_{b-1}$ and $N_b$ are as defined in Eq.~(\ref{nbnb1}).

As discussed in Sec.~\ref{s:derive}, the MSD boundary condition can be expanded out, expressing the $b-1$ Laplacian in terms of the internal scheme in order to be able to evaluate the MSD boundary condition simultaneously with the interior points (when using explicit time-stepping schemes, this is usually unnecessary).  As an example, using central-differencing in space for the one-dimensional NLSE, Eq.~(\ref{MSDlap1}) on a uniform grid can be expanded as
\begin{equation}
\nabla^2\Psi_b \approx 
\left[ 
\frac{1}{h^2}
\left(
\mbox{Im}\left(i\, \frac{\Psi_b + \Psi_{b-2}}{\Psi_{b-1}}\right) - 2
\right)
+ \frac{1}{a}\,\left(N_{b-1} - N_b\right)\right]\Psi_b,
\end{equation}
where $N_{b-1}$ and $N_b$ are as defined in Eq.~(\ref{nbnb1}) and $h$ is the cell spacing.

\section{Numerical Results}
\label{s:num}
In order to demonstrate the usefulness and advantages of the MSD boundary condition, we show a few example simulations of the NLSE.   The MSD boundary condition is compared to a Laplacian-zero (L0) boundary condition (defined as $\nabla^2 \Psi_b = 0$) as it is an easy-to-implement boundary condition that may be used for constant background-density simulations.

Every numerical boundary condition, including the MSD, is only valid under its governing assumptions.  Therefore, if the modulus-squared of a solution is changing at the boundary, the MSD is not expected to work well, just as if the Laplacian of the wavefunction is far from zero at the boundary, the Laplacian-zero boundary condition is not expected to do well.  Therefore, comparisons of which boundary condition is best is very often problem-specific.  That being said, comparing the MSD to the L0 is still valuable in that for some problems, both boundary conditions are suitable, allowing for a fair comparison.  In addition, since limiting the size of the required grid is very important (especially for higher-dimensional simulations), it is useful to see which boundary condition allows for the use of the smallest (tightest) grid within acceptable accuracy limits.

While there are numerous numerical methods that can be used to simulate solutions to the NLSE, we perform the simulations using the code package NLSEmagic\footnote{Available at \url{http://www.nlsemagic.com}} \cite{NLSEmagic} which uses the fourth-order in time Runge-Kutta method with either a second-order central difference (RK4+CD) or a two-step high-order compact fourth-order central-differencing in space (RK4+2SHOC) \cite{ME_2SHOC} (all simulations in this paper are performed using the RK4+2SHOC scheme).

\subsection{One-Dimensional Dark Solitons in the NLSE}
\label{s:num1d}
For the numerical one-dimensional tests we use a co-moving dark soliton solution with $V({\bf r})=0$ given by \cite{BEC_RCbook}:
\begin{equation}
\label{1dmds}
\Psi(x,t) = \sqrt{\frac{\Omega}{s}}\,\mbox{tanh}\left[\sqrt{\frac{-\Omega}{2a}}(x-ct)\right]\mbox{exp}\left(i\left[\frac{c}{2a}x + \left(\Omega - \frac{c^2}{4a}\right)t\right]\right)
\end{equation}
where $c$ is the velocity of the soliton and $\Omega$ is a chosen parameter representing the frequency of the solution.

The first test is to compare the error for a stationary dark soliton (setting the co-moving velocity $c$ to $0$).  In such a case, for a large enough domain, the L0 boundary condition can be used since $\Psi$ flattens out as $|x|\rightarrow\infty$.  In contrast, the MSD boundary condition should work on any sized domain, since its underlying assumption of constant density is valid anywhere along the steady-state solution.  A one-sided second-order differencing (1SD) boundary condition defined in one-dimension as
\[
\frac{\partial^2 \Psi}{\partial x^2} \approx \frac{1}{h^2} \, \left(-\Psi_{b-3} + 4\,\Psi_{b-2} -5\,\Psi_{b-1} + 2\,\Psi_b\right),
\]
is also used for comparison in this case.

\begin{figure}[htb!]
\centering
\includegraphics[width=2.3in]{ 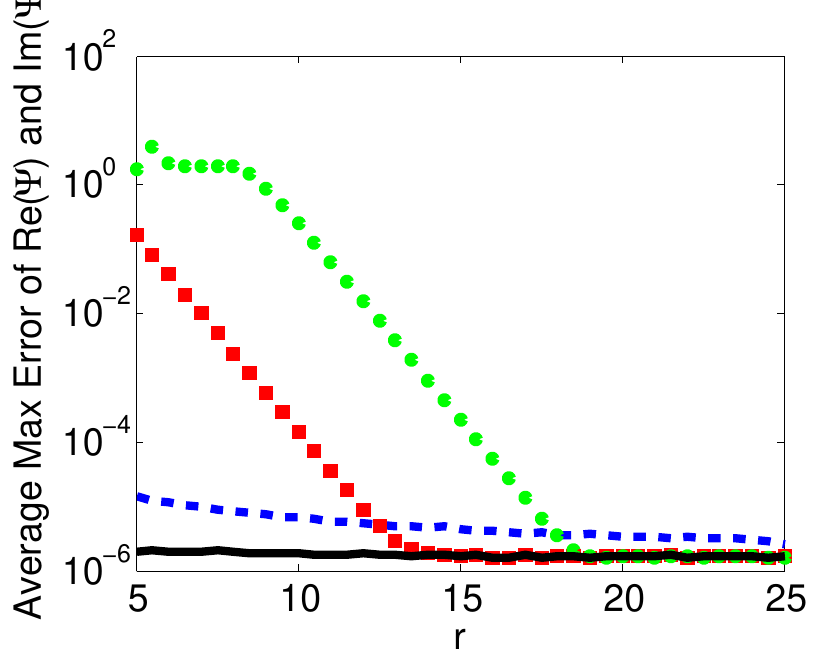}
\includegraphics[width=2.3in]{ 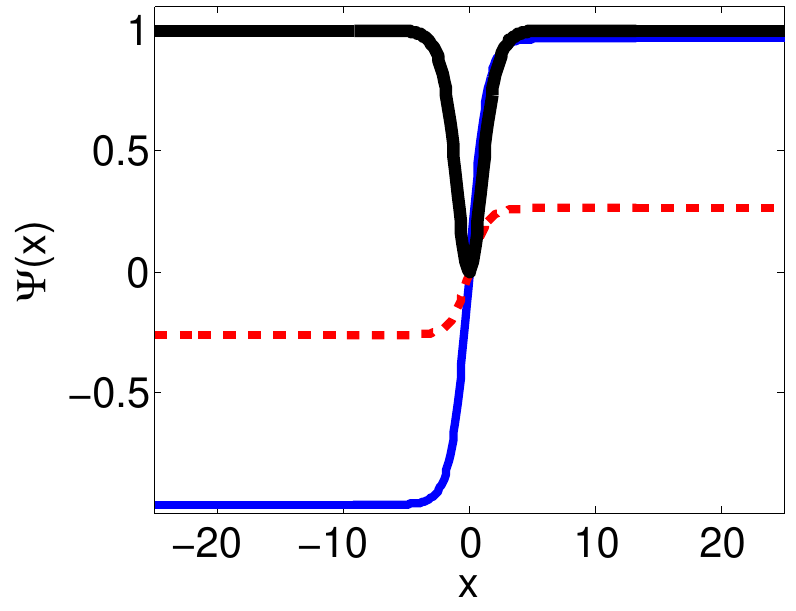}
\caption[MSD boundary condition error comparison for a one-dimensional dark soliton]{(Color online) Left:  Comparisons of errors for simulating the dark-soliton solution (with $c=0$, $\Omega=-1$, $a=1$, and $s=-1$) of Eq.~(\ref{1dmds}) using MSD (blue dashed line), L0 (red squares), 1SD (green dots), and exact (black solid line) boundary conditions for $r\in[5,5.5,6,...,25]$.  The simulations are run to an end time of $t=50$ with a time-step of $k=0.005$ and a spatial step of $h=0.1$.  Right:  Depiction of the soliton for the maximum domain ($r=25$) simulation at time $t=50$.  The thin (blue) and dashed (red) lines are the real and imaginary parts of $\Psi$ respectively, while the thick solid (black) line is the modulus-squared, $|\Psi|^2$. \label{f:1dc0}}
\end{figure}

We define a radius, $r$, which represents the distance from the center of the soliton to the edge of the computational domain.  Simulations are performed for various lengths of $r$ ranging from $r=5$ (approximately equal to the width of the soliton) to $r=25$ (the distance where the boundary value of the modulus of the soliton is approximately equal to the infinite background density ($\Omega/s$) minus machine epsilon $\epsilon \approx 10^{-16}$).  The average of the maximum errors in the real and imaginary parts of $\Psi$ over the length of the simulations are recorded.  The results are shown in Fig.~\ref{f:1dc0}.
Keeping in mind that $h^4=10^{-4}$ and that the spatial scheme has an accuracy of $O(h^4)$, it is clear that the MSD boundary condition performs well even when the domain is small.  The L0 and 1SD have less error for large $r$ and converge to the error from the exact boundary conditions.  This is understandable since as the Laplacian tends towards zero rapidly as $r$ increases, the L0 and 1SD both equate the Laplacian to $0$ and have no additional errors associated with them.  Even so, the MSD can simulate the solution to acceptable accuracy at a much smaller grid size than the L0 or the 1SD, demonstrating its usefulness in this case.

\begin{figure}[htb!]
\centering
\includegraphics[width=2.5in]{ 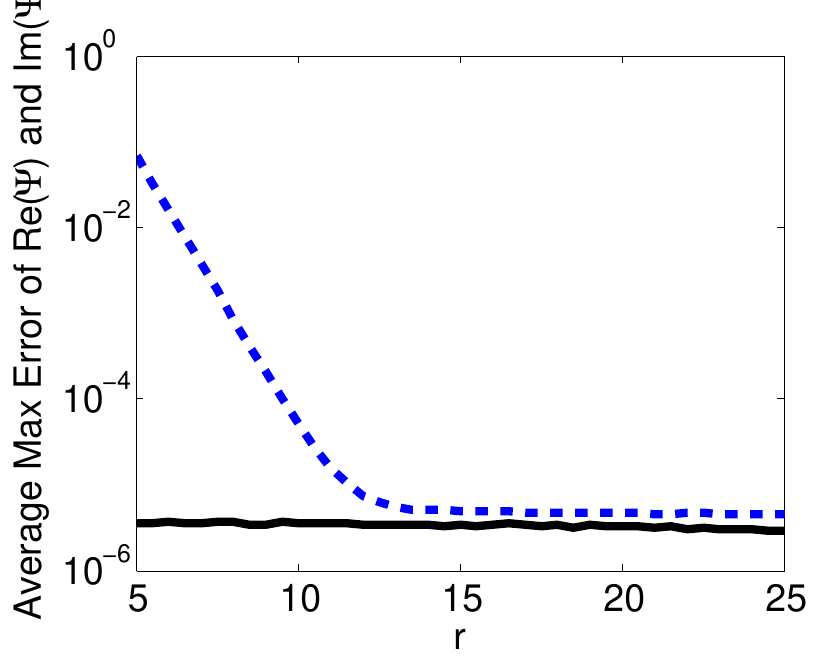} \\
\includegraphics[width=1.25in]{ 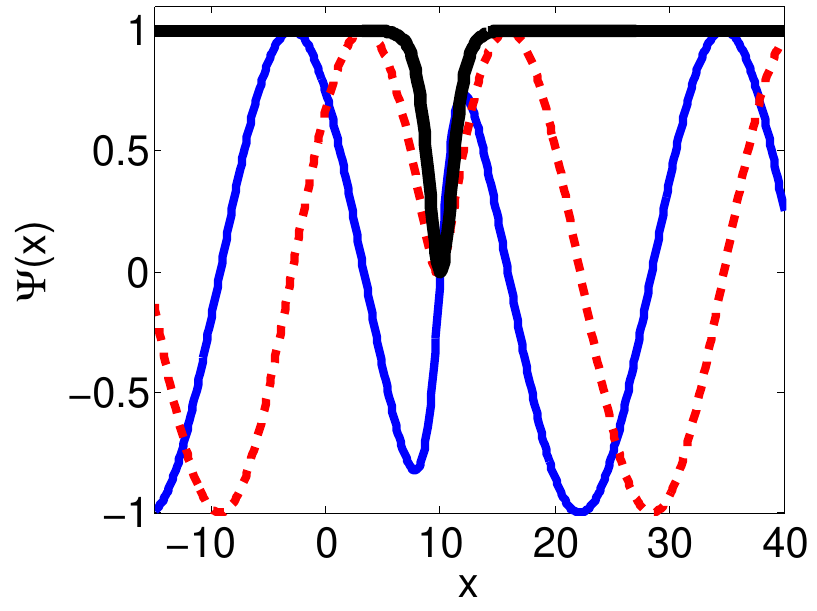}
\includegraphics[width=1.25in]{ 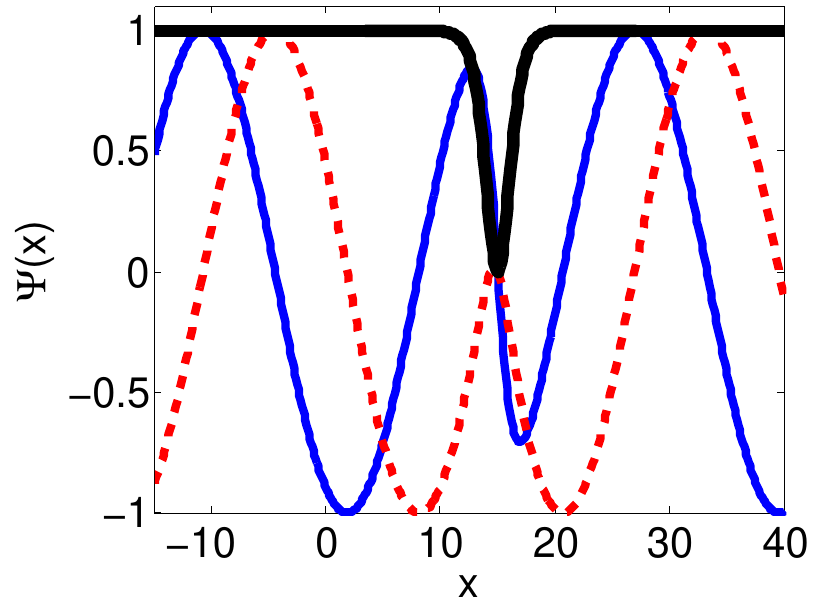}
\includegraphics[width=1.25in]{ 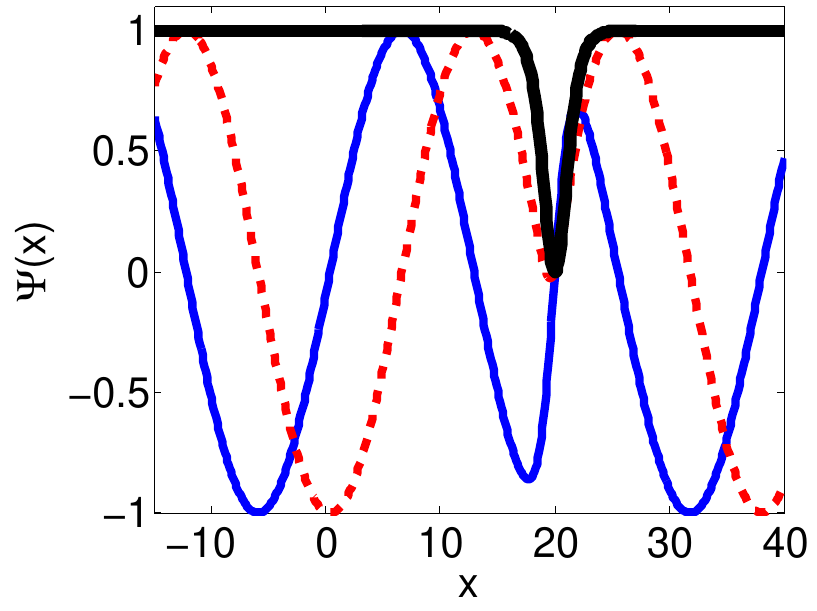}
\includegraphics[width=1.25in]{ 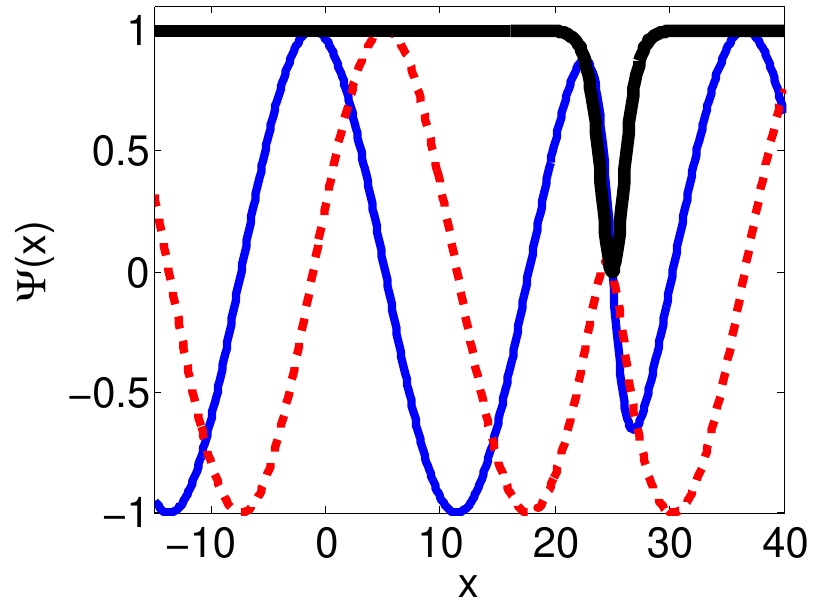}
\caption[MSD boundary condition error comparison for a one-dimensional co-moving dark soliton]{(Color online) Top:  Comparisons of errors for simulating the dark-soliton solution of Eq.~(\ref{1dmds}) using MSD (blue dashed line), and exact (black solid line) boundary conditions for various domain sizes.  The simulation parameters and figure descriptions are the same as in Fig.~\ref{f:1dc0} except that here, the velocity is $c=0.5$.  Bottom:  Depiction of the soliton during the simulation with $r=15$ at times $t=20, 30, 40, 50$ using the MSD boundary condition.  \label{f:1dc05}}
\end{figure}

For the next test, the soliton is given a co-moving velocity of $c=0.5$.  Since the L0 assumptions are completely invalid at any domain size, it is not used for comparison (the 1SD is also not used as it fails quickly at any domain size as well).  Therefore, the MSD boundary condition is only compared to the exact boundary condition.
In this case the domain is set to be a distance $r$ to the left of the initial position of the soliton as before, and a distance $r+c\,T$ to the right (where $T$ is the simulation end-time) in order to account for its movement.  In Fig.~\ref{f:1dc05}, we show the results of the simulations.
We see that the MSD boundary condition performs extremely well as long as the soliton is far enough from the boundaries (in this case `far enough' is about equal to where the background density minus $|\Psi|^2$ at the boundary is approximately $h^4$).  Through our review of the literature, it seems that the MSD is the \emph{only} simple-to-implement boundary condition that can handle such a co-moving back-flow.  We note that in other tests not presented here using the RK4+CD scheme, it was found that the MSD can have \emph{less} error than the exact boundary conditions for large values of $r$.  This observation can possibly be explained by the fact that using exact boundary conditions in fourth-order Runge-Kutta schemes can actually introduce errors in the solution as described in Ref.~\cite{FD_PS_BOOK}.

\subsection{Two-Dimensional Dark Vortices in the NLSE}
\label{s:num2d}
To test the MSD boundary condition in a more complicated setting, we use the known dynamics of dark vortices in the two-dimensional NLSE.  Dark vortices are described as
\begin{equation}
\label{2dvortex}
\Psi(r,\theta,t) = f(r)\,\exp\left[i\,(m\,\theta + \Omega \,t)\right]
\end{equation}
where $m$ is the vortex charge (also known as the winding number) and $\Omega$ is the frequency which is directly related to the background density, $\rho$, as $\rho=\Omega/s$.  The real-valued radial profile $f(r)$ can be obtained numerically by inserting Eq.~(\ref{2dvortex}) into Eq.~(\ref{NLSE}) and solving the resulting ODE for $f(r)$ using a nonlinear equation solver (in our case, a Newton-Krylov GMRES($m$) solver in a package called {\tt nsoli} \cite{OPT_NSOLI_BOOK}).  As an initial iterate for the solver, we use the asymptotic profile approximation given by \cite{BEC_DYN_NONLIN}
\begin{equation}
\label{vprofile}
f(r>0)\approx \mbox{Re}\left[\sqrt{\frac{\Omega}{s} + \frac{a\,m^2}{s\,r^2}}\right].
\end{equation}
It is important to note that the tails of the dark vortices converge to the background density much slower than the one-dimensional dark solitons of Eq.~(\ref{1dmds}).  For example, in the one-dimensional dark soliton, we could extend the domain to a size so that the wave-function is at a value of $\sqrt{\Psi_\infty}-\epsilon$ ($\epsilon\approx 10^{-16}$), and this would give a radius of around $25$ (for our parameter choices).  To get the same boundary value in a two-dimensional vortex of charge $m=1$, we would require a radius of over \emph{1 million}!  Therefore, the ability of the chosen boundary condition to not effect the dynamics of the system on a small grid is vital.

\begin{figure}[htb!]
\centering
\includegraphics[width=1.72in]{ 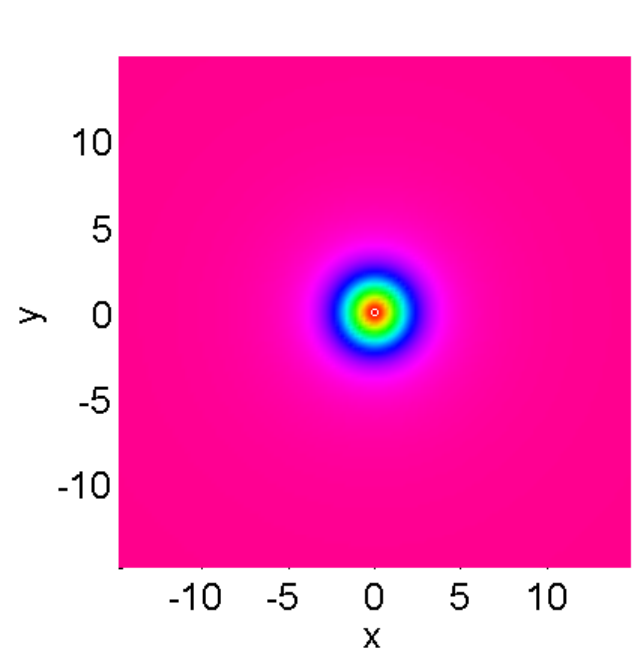}
\includegraphics[width=1.72in]{ 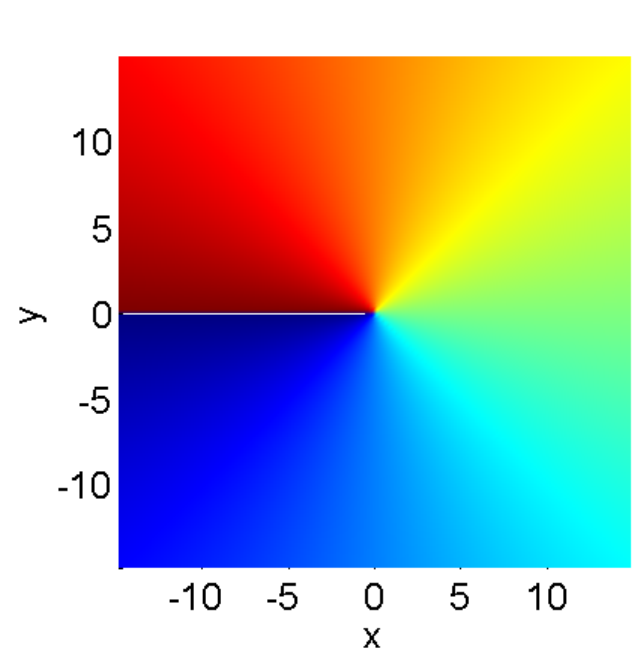}
\\
\includegraphics[width=1.25in]{ 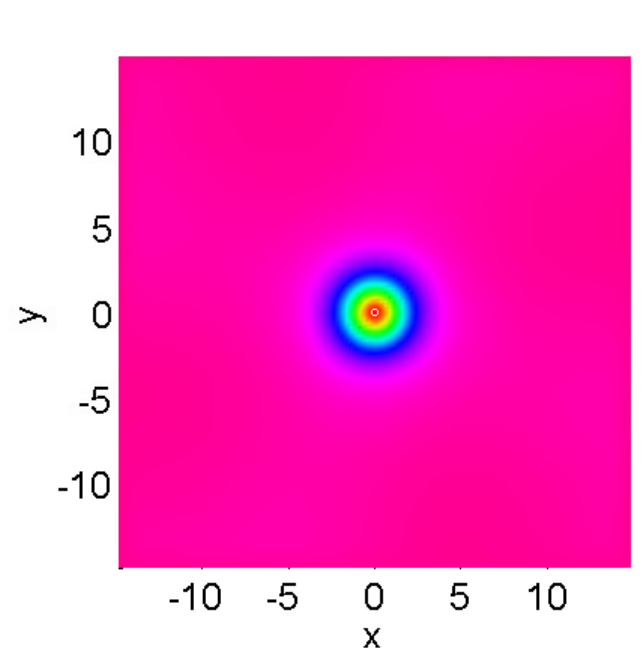}
\includegraphics[width=1.25in]{ 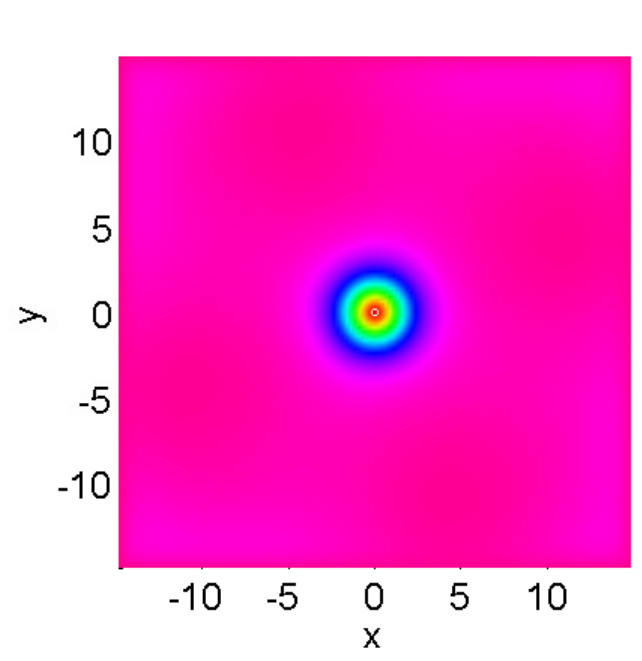}
\includegraphics[width=1.25in]{ 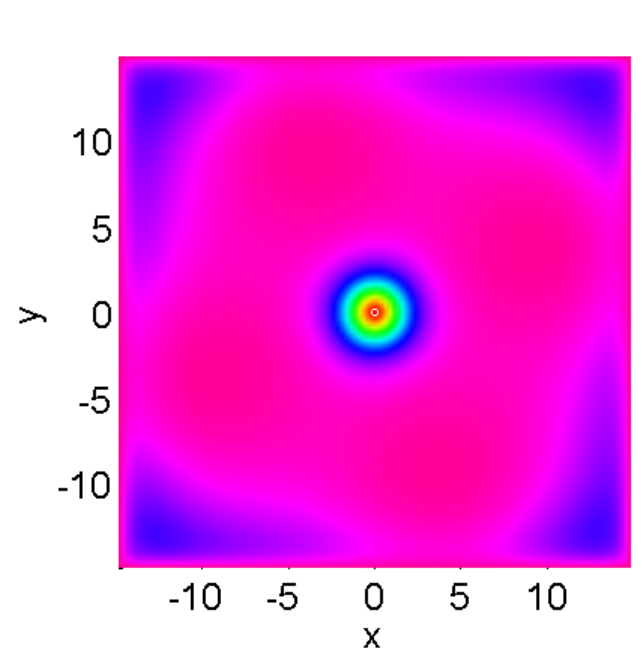}
\includegraphics[width=1.25in]{ 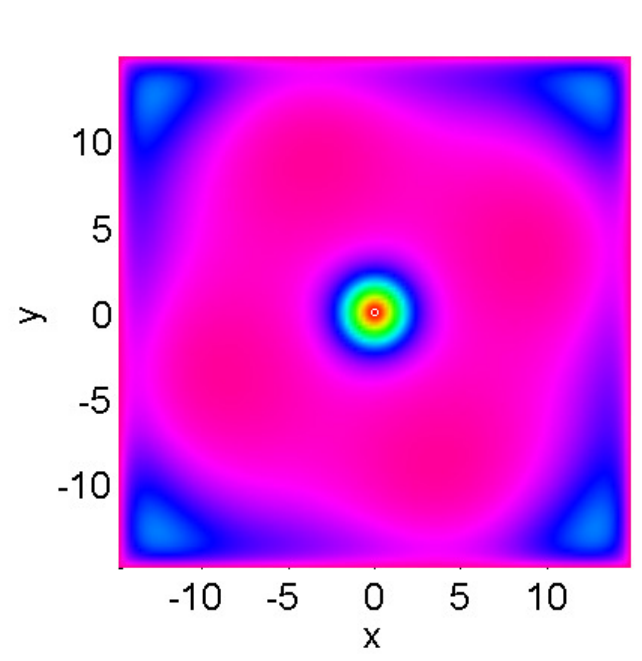}
\\
\includegraphics[width=1.25in]{ 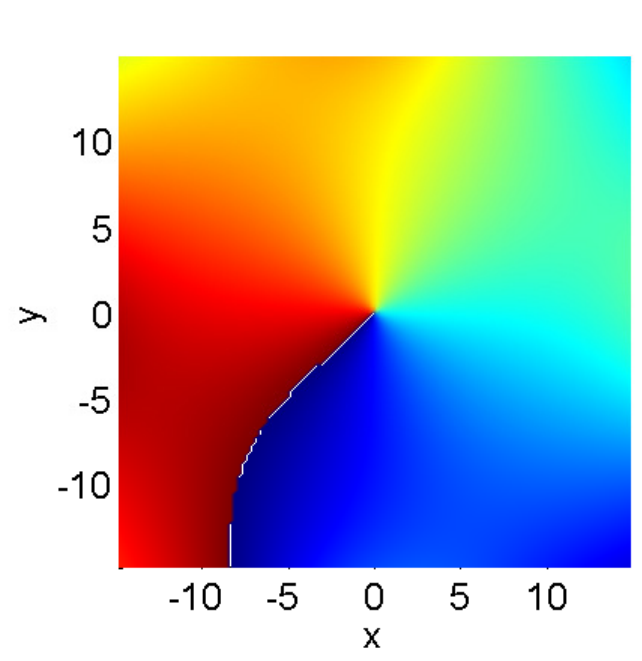}
\includegraphics[width=1.25in]{ 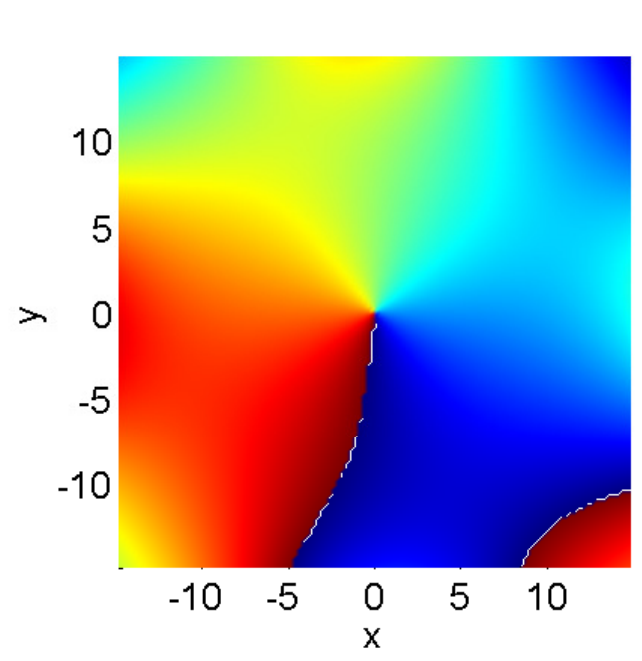}
\includegraphics[width=1.25in]{ 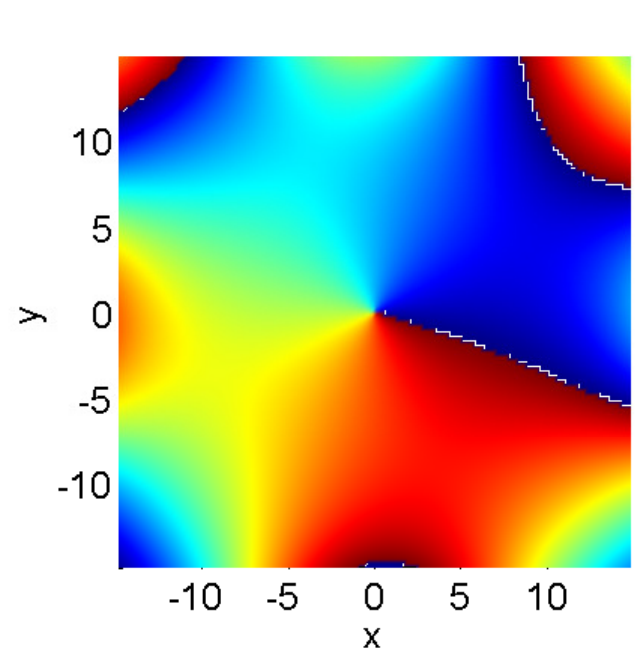}
\includegraphics[width=1.25in]{ 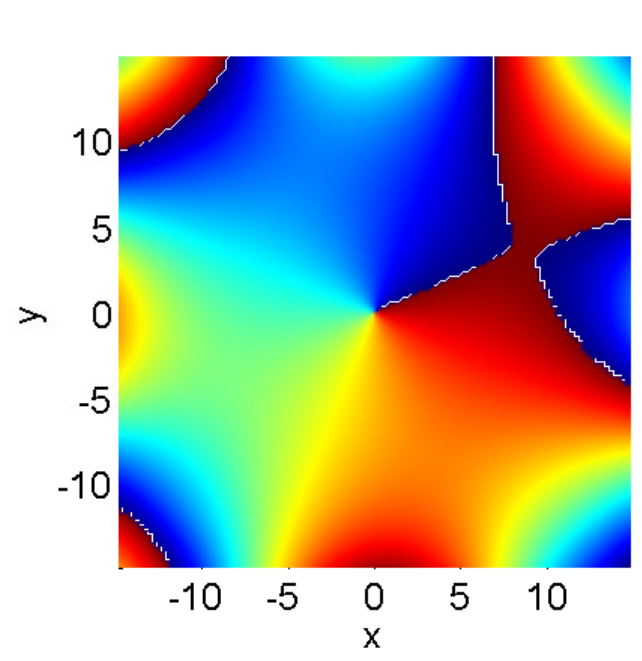}
\\
\includegraphics[width=1.72in]{ 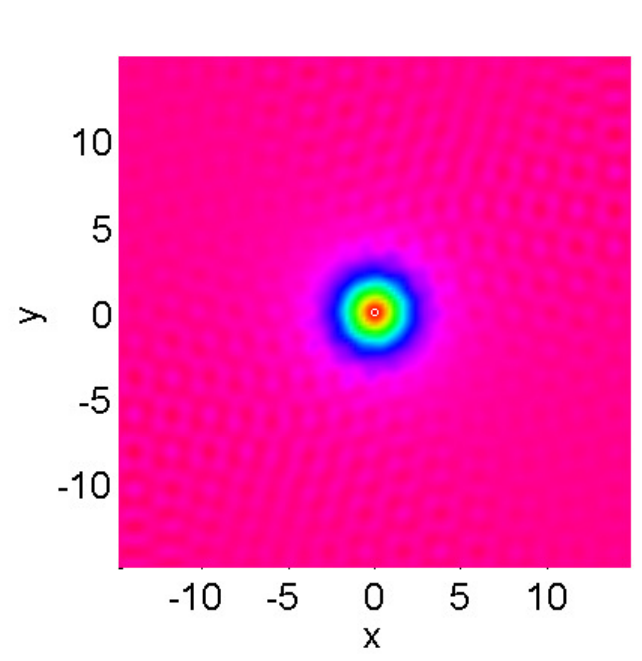}
\includegraphics[width=1.72in]{ 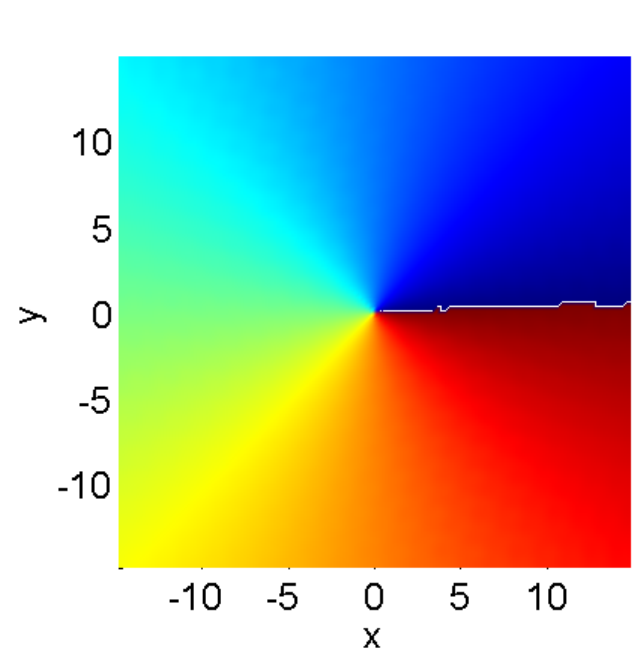}
\caption[MSD boundary condition dark vortex phase test.]{(Color online) Top row: Modulus-squared and phase of the initial condition of a single dark vortex solution.  Second (third) row: Snapshots of the density (phase) using the Laplacian-zero boundary condition at times $600$, $1200$, $2400$, and $3000$.  Bottom row: Modulus-squared and phase of simulation using the MSD boundary condition at time $t=50,000$.  All simulations use a spatial-step of $h=0.25$ and a time-step of $k=0.001$ on a grid size of $120 \times 120$.  \label{f:2d1vort}}
\end{figure}

The first test is to simulate a single unitary-charged ($m=1$) vortex which is known to be a stable steady-state solution to the NLSE \cite{BEC_RCbook}.  Choosing a moderate domain size (a $120 \times 120$ grid with a spatial step-size of $h=0.25$), we integrate the NLSE for considerable time (up to $T=50,000$) using both the MSD and L0 boundary conditions.  The results are shown in Fig.~\ref{f:2d1vort}.
It is clear that the L0 boundary condition is only useful for shorter simulations, as it quickly suffers from phase discrepancies which get worse as time progresses, until the point where the solution is distorted badly enough to be unusable.  From successive tests it is found that this effect occurs for the L0 boundary condition for even very large domains, but the time of the onset of the distortions is delayed longer as the domain size is increased.  In contrast to this, the MSD boundary condition creates {\em no distortions} in the phase or modulus-squared of the solution for very long simulations (in this case up to $t=50,000$).

\begin{figure}[htb]
\centering
\includegraphics[width=2.5in]{ 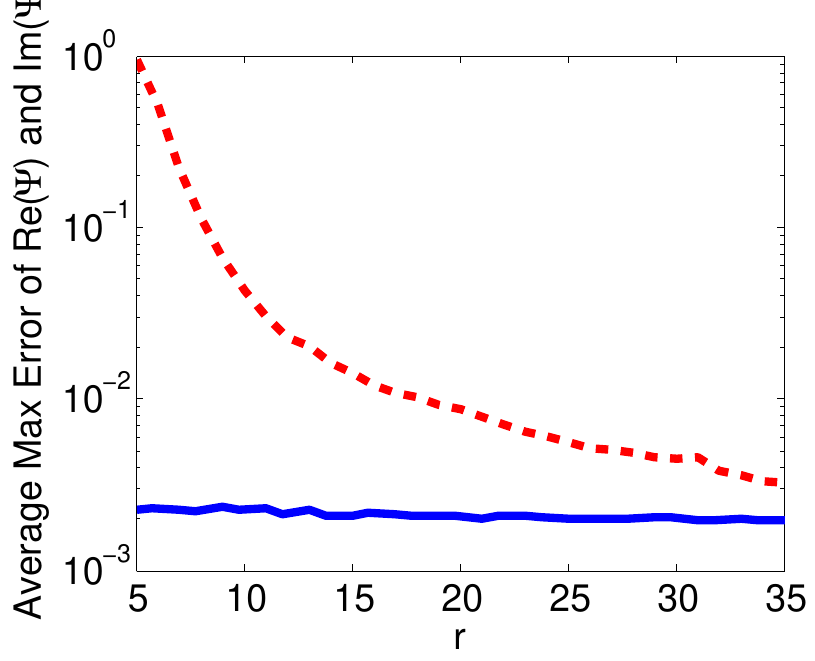}
\caption[MSD boundary condition dark vortex grid requirements test]{(Color online) Maximum error of the modulus-squared of a steady-state dark vortex using the L0 (red dashed line) and the MSD (solid blue line) boundary conditions for various domain sizes.  The initial numerically optimized solution of the vortex is used as the `true solution' for error comparisons.  The simulations were run to an end-time of $t=300$ with a spatial-step of $h=0.25$ and time-step $k=0.01$.  \label{f:2d1vort_r}}
\end{figure}

Given a set end-time, it is useful to determine the size of the grid needed to properly simulate the vortex for a given boundary condition.  Using the initial vortex solution as representing the `exact' steady-state solution, we can track the error in the modulus-squared of the solution over the course of the simulations as the grid-size is increased.  A radius $r$ is defined as the distance from the center of the vortex to the edge of the domain in the $x$ or $y$ direction (whichever is smallest).  Figure~\ref{f:2d1vort_r} shows the results of varying $r$ from $5$ to $35$ for a simulation with an end-time of $t=300$.
It is clearly seen that a much larger grid is required for the L0 boundary condition to be close to the effectiveness of the MSD boundary condition in this case.  For example, even at a large $280 \times 280$ grid, the L0 boundary condition did not have as low an error as the MSD boundary condition did on a $41 \times 41$ grid!  Once again, this is understandable considering that the profile of the vortex does not flatten out rapidly as in the one-dimensional case.

\begin{figure}[htb!]
\centering
\hskip-0.25cm
\includegraphics[width=1.5in]{ 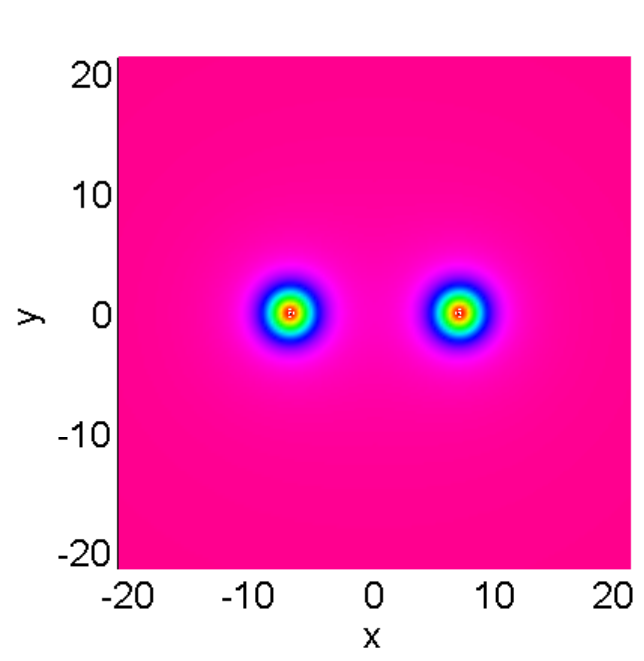}
\hskip0.25cm
\includegraphics[width=1.5in]{ 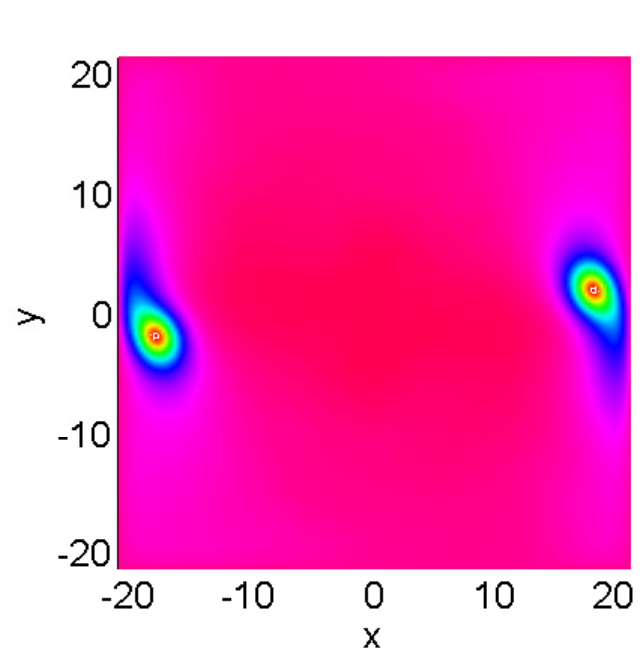}
\hskip0.25cm
\includegraphics[width=1.5in]{ 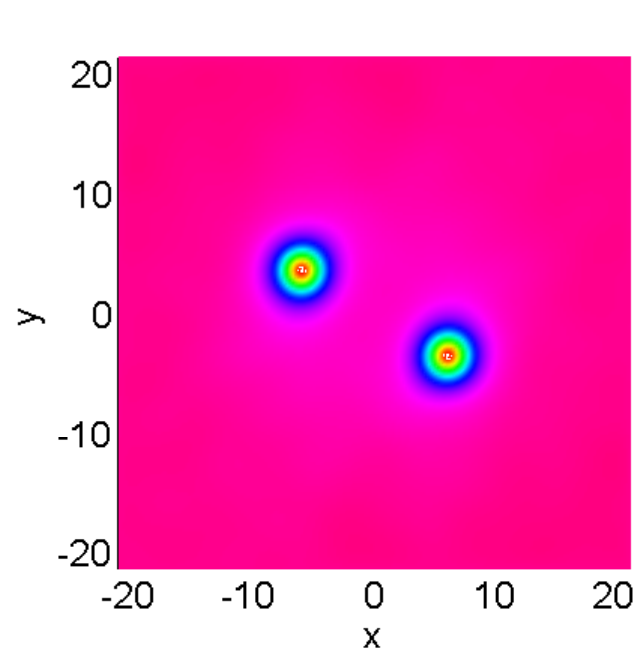}
\hskip0.25cm
\\
\includegraphics[width=1.4in]{ 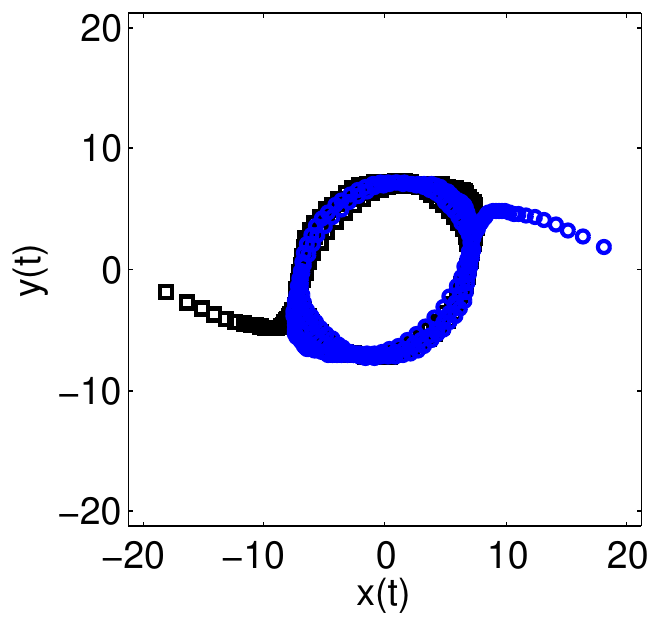}
\includegraphics[width=1.7in]{ 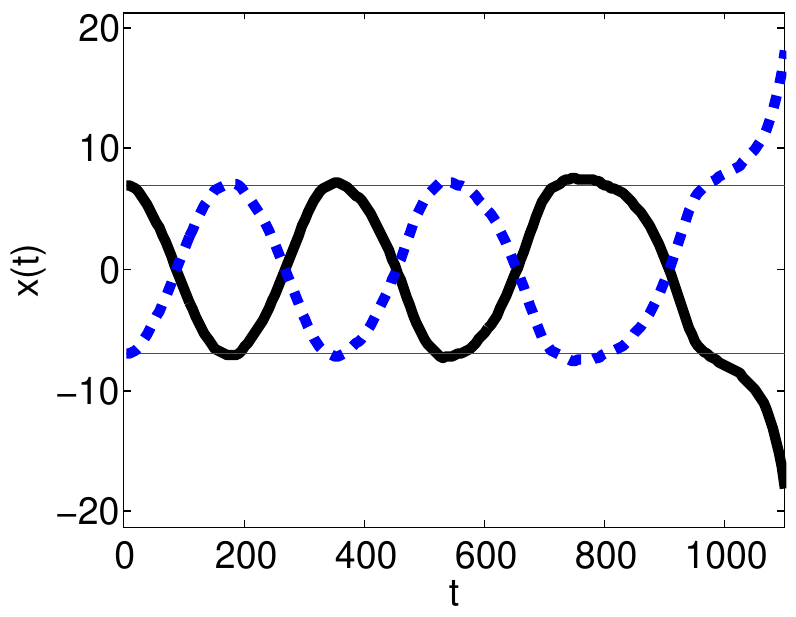}
\includegraphics[width=1.7in]{ 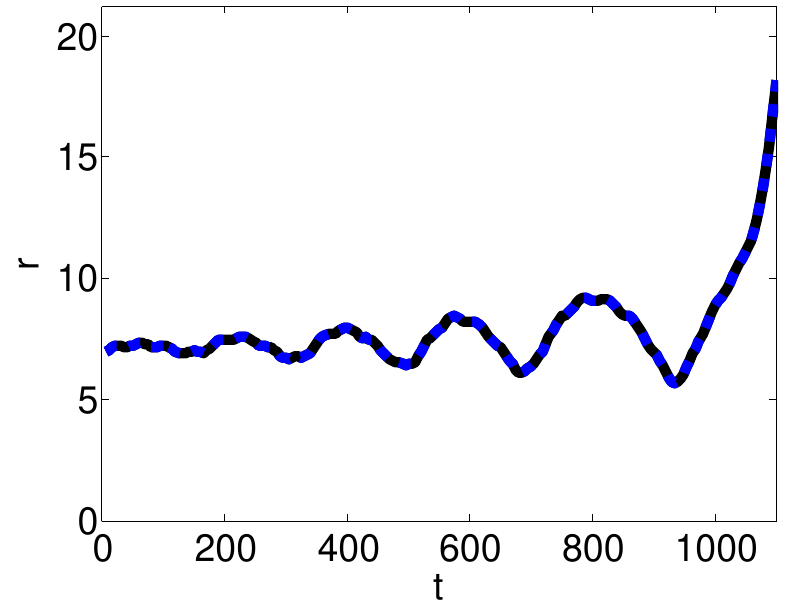}
\\
\includegraphics[width=1.4in]{ 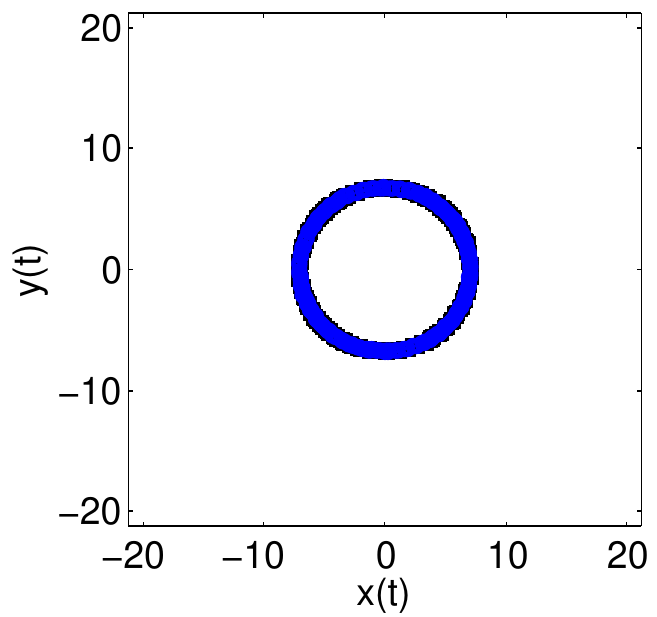}
\includegraphics[width=1.7in]{ 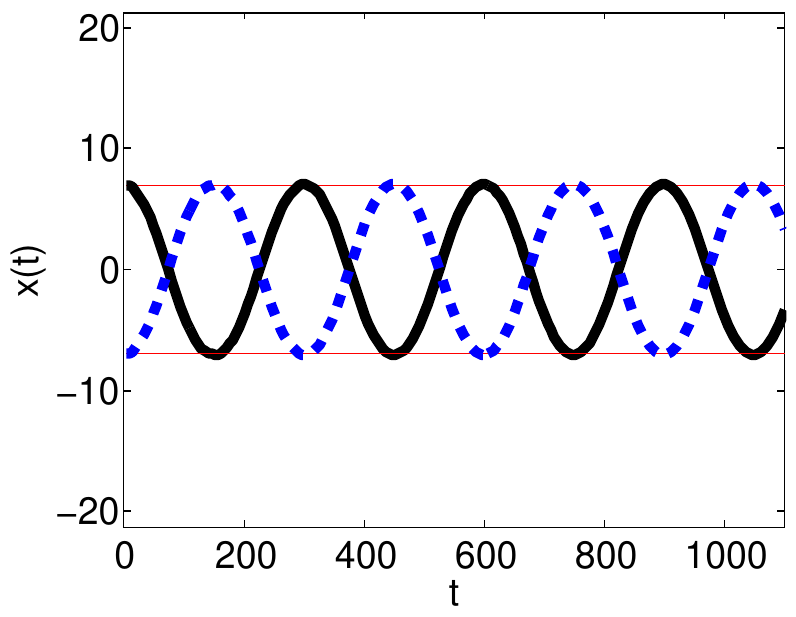}
\includegraphics[width=1.7in]{ 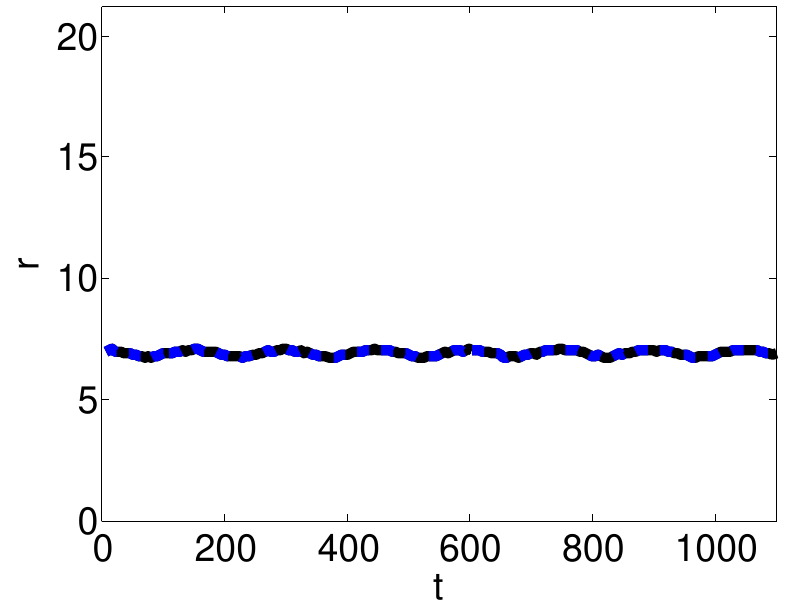}
\\
\hskip0.275cm
\includegraphics[width=1.4in]{ 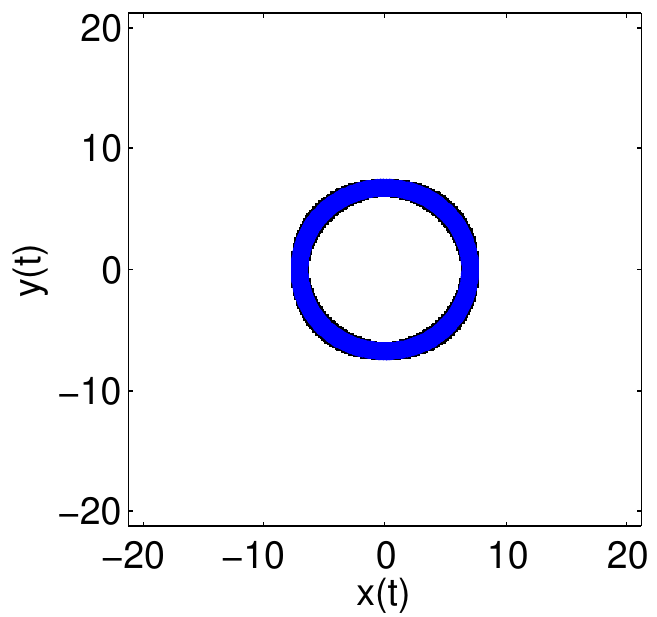}
\includegraphics[width=1.81in]{ 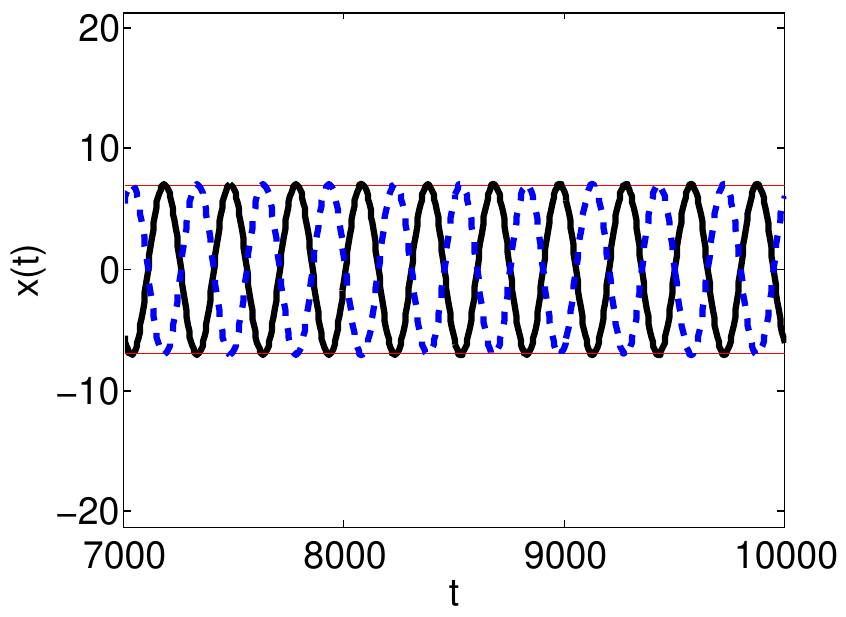}
\hskip-0.23cm
\includegraphics[width=1.81in]{ 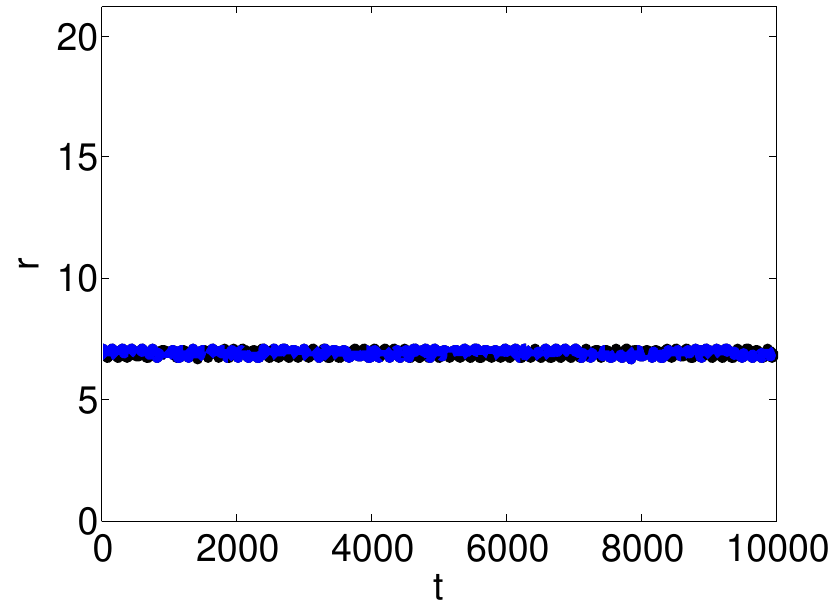}
\caption[MSD boundary condition rotating dark vortices test.]{(Color online) Top row, left to right: Modulus-squared of initial condition of two single-charge vortices separated by a distance of $7$ from the center of the grid, a snapshot of the simulation using the L0 boundary condition at time $t=1100$, and a snapshot using the MSD boundary condition at time $t=10,000$.  Second row, left to right: Traced positions of the two vortices (square and circles respectively) over the course of a simulation with an end time of $t=1100$ using L0 boundary conditions, $x$ positions of the two vortices (solid and dashed line respectively) versus time, and the computed distance between the vortices (solid and dashed line respectively) and the center of the grid for an end time of $t=1100$. Third row:  Same description as the second row, but using the MSD boundary condition.  Bottom row:  Same as in the third row, but with an end-time of $t=10,000$ (only times $7,000$ through $10,000$ shown in the $x$ versus $t$ plot).  In all simulations the spatial-step is $h=0.25$ and the time-step is $k=0.01$ with a grid size of $171 \times 171$.\label{f:2d2vort}}
\end{figure}

As an additional example to test the MSD boundary condition, we simulate two equal-charge vortices whose interaction is known to produce a rotating circular motion of the two vortices orbiting each other \cite{Bao_SJAM_07,Bao_JCP_13,Bao_EJAM_07,VORT_DYN_REC_DEV}.  Using a fixed grid size of $171 \times 171$, the simulations are run for long times using the L0 and MSD boundary conditions.  The results are shown in Fig.~\ref{f:2d2vort}.
We see that once again, using the L0 boundary condition causes a break-down in the dynamics, eventually causing the two vortices to decouple from each other and fling into the boundaries.
In contrast, the MSD boundary condition allows for near-perfect rotational dynamics for indefinite simulation times, even for small grid sizes.  It is also noticeable that the period of rotation of the vortices is shorter when using the MSD boundary condition when compared to using the L0 boundary condition.  Through further simulations with larger grid sizes, we have observed that the period of rotation converges to the same value (between the two values of the period displayed) for both boundary conditions at approximately the same grid-size (around $250 \times 250$).  Therefore, the MSD boundary condition performs as well as the L0 in terms of rotational period, but much better in terms of long-term dynamics and minimum grid-size requirements.

\begin{figure}[htb!]
\centering
{\includegraphics[width=3in]{ 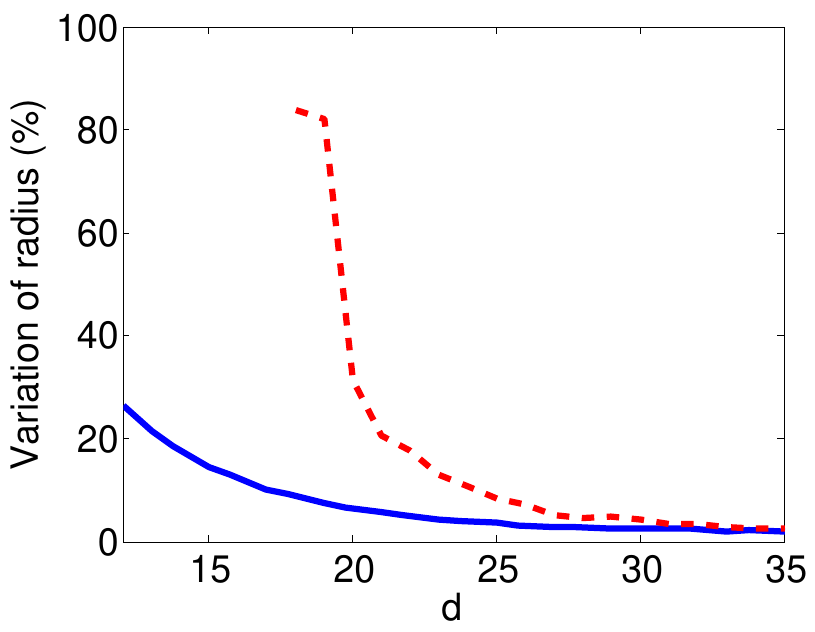}}
\\
{\includegraphics[width=2in]{ 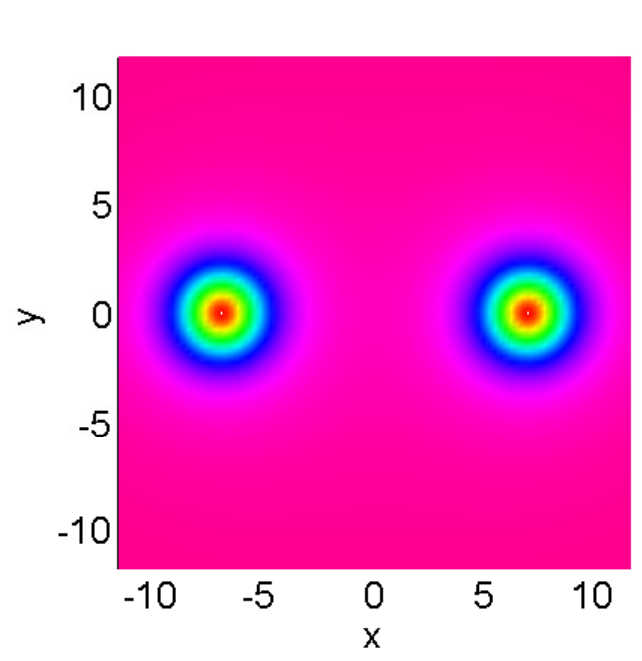}}
{\includegraphics[width=2in]{ 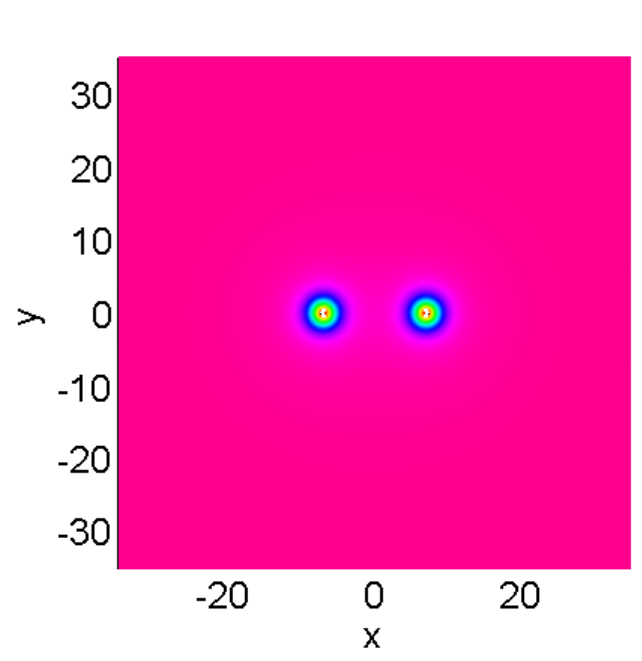}}
\caption[MSD boundary condition rotating dark vortices test.]{(Color online) Top:  Percentage of radius variation when compared to the initial radius of the vortices to the center of the grid as a function of grid size $d$ (defined as the distance from the center of the grid to its edge along the $x$ or $y$ direction) for a simulation with an end-time of $t=480$ using both the L0 (dashed (red) line) and the MSD (solid (blue) line) boundary conditions.  The grid size, $d$, is varied from $12$ (the minimum size required to resolve both vortices adequately) to $35$.  Results for the L0 boundary condition below $d=18$ are not shown as the vortices hit the grid wall before the simulation ends.  The other simulation parameters are the same as those in Fig.~\ref{f:2d2vort}.  Bottom:  Modulus-squared of the initial condition for grid-sizes $d=12$ and $d=35$. \label{f:2d2vort2}}
\end{figure}

If the end-time of the simulation is fixed to be such that the vortices will rotate at least one complete rotation (here we use $t=480$), we can record the deviation from perfect circular motion as the grid size (given as the distance, $d$, from the center of the grid to the edge along the $x$ or $y$ direction) is varied.  The vortices are tracked during the simulations and the maximum variation of radius compared to a constant radius to the center of the rotating vortices is recorded.  The results are shown in Fig.~\ref{f:2d2vort2}.
We see that the L0 boundary condition requires a very large grid size to capture the correct dynamics (and completely fails for smaller grids), while the MSD is able to capture the dynamics to an acceptable degree on a much smaller grid-size.  The discrepancies when using the MSD boundary condition at low grid sizes (up to $20\%$ radius variation) is understandable since at those distances the boundaries become far from steady-state due to the perturbations in the background density caused by the motion of the vortices.

\subsection{Three-Dimensional Dark Vortex Rings in the NLSE}
To further show the usefulness of the MSD boundary condition in
multi-dimensional settings, we compare the MSD boundary condition
to the L0 within simulations of three-dimensional dark vortex rings in the NLSE.
Following the natural progression of complexity, starting from dark
solitons in 1D and vortices in 2D, let us now turn our attention
to their 3D analogs: vortex rings (VRs).
A vortex line is the 3D extension of a 2D vortex by infinitely
extending the solution into the axis perpendicular to the vortex plane.
%
%
If a vortex line is bent enough to close on itself or if two
vortex lines are close enough to each other they can produce
a vortex ring \cite{VR:CrowInstab}.
VRs are 3D structures whose core is a closed loop with vorticity
around it \cite{VR_HEL_BOOK} (i.e., a vortex that is looped back into itself).
%
%

%
Numerical studies directly focusing on the structure and stability of dark vortex rings in trapped BECs have been addressed.  These include the generation of steady-state vortex rings in toroidal \cite{VR_BEC_TOROIDAL} and harmonic \cite{VR_BEC_STRUCTGOOD, VR_BEC_LINESOURCE} trapping potentials, multiple parallel ring structures \cite{VR_BEC_PARALLEL}, as well as the study of the bending-wave instability of vortex rings \cite{VR_BEC_BEND_STB} and their degeneration into quantum turbulence \cite{VR_BEC_BEND_TURB}.  Some dynamics of vortex rings in trapped BECs have been explored, such as their oscillatory motion resulting in self-annihilation \cite{VR_BEC_VEL_DYN}, and collisions of vortex rings with solitons in cylindrical BECs \cite{VR_BEC_SOLSMASHTUBE, VR_BEC_TUBE_COLLIDEEU}.

From a numerical perspective, a salient feature of VRs is that they posses
an intrinsic transverse velocity associated with them \cite{VR_NLSE_VEL_71}. Therefore, in order to run long simulations of the rings (for instance, to study stability under small perturbations, or interactions between multiple co-moving rings), a very large grid size is required.  Often, due to the size of the simulations, a large enough grid is not within the memory limitations of the computers being used for the simulations.  To avoid this problem, a vortex ring can be made to be a `steady-state' by applying a background velocity equal and opposite to the vortex ring's intrinsic velocity.  By doing this, long-term simulations of the rings can be performed, but with many fewer grid points.

\begin{figure}[htb!]
\centering
\includegraphics[width=4.5in]{ 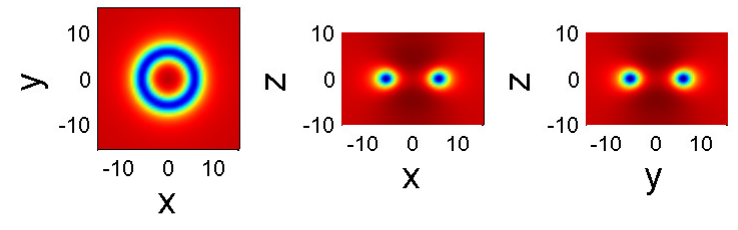}
\\
\includegraphics[width=4.5in]{ 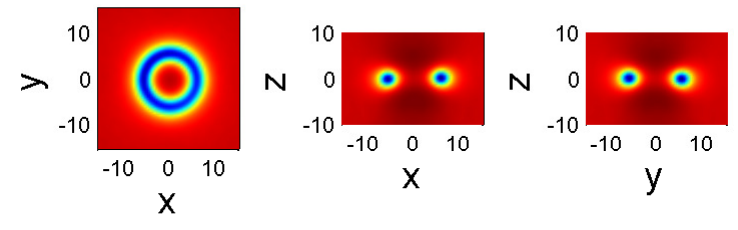}
\caption[Simulation of a steady-state vortex ring using the MSD boundary condition.]{(Color online) Two-dimensional cuts of the modulus-squared of the initial condition of a steady-state vortex ring of radius $6$ amidst a back-flow.  Bottom row:  Two-dimensional cuts of the same vortex ring at simulation time $t=1500$.  The simulation uses a spatial step-size of $h=0.5$, and a time-step of $k = 0.035$. \label{f:vr}}
\end{figure}

In Fig.~\ref{f:vr} we show a simulation of a steady-state vortex ring amidst a back-flow at various time intervals using the MSD boundary condition.  The vortex ring solution is found by seeding the numerically-exact two-dimensional vortex described in Sec.~\ref{s:num2d} (with the correct back-flow velocity added to it) into a nonlinear optimization routine utilizing the two-dimensional axisymmetric version of the three-dimensional NLSE (for details, see Ref.~\cite{ME_DISS}).
As was the case with the co-moving dark soliton of Fig.~\ref{f:1dc05}, the assumptions underlying the L0 boundary condition are not valid at any grid size, and so it is not tested in this specific case.
The vortex ring remains stationary for very long simulation times (up to $t=1500$ in this case), further demonstrating the MSD boundary condition's usefulness in studying co-moving solutions.

As was the case for two-dimensional vortices, when running dynamical simulations of vortex rings, it is important to know how large to set the computational grid in order to minimize the effect of the boundary conditions on the interior dynamics.  Since the vortex rings have a constant transverse velocity, a good test of the effects of the boundary conditions is to track a vortex ring in a simulation and observe the velocity for various grid sizes.  As the size of the grid is expanded, this velocity should converge.

\begin{figure}[htb!]
\centering
\includegraphics[width=1.65in]{ 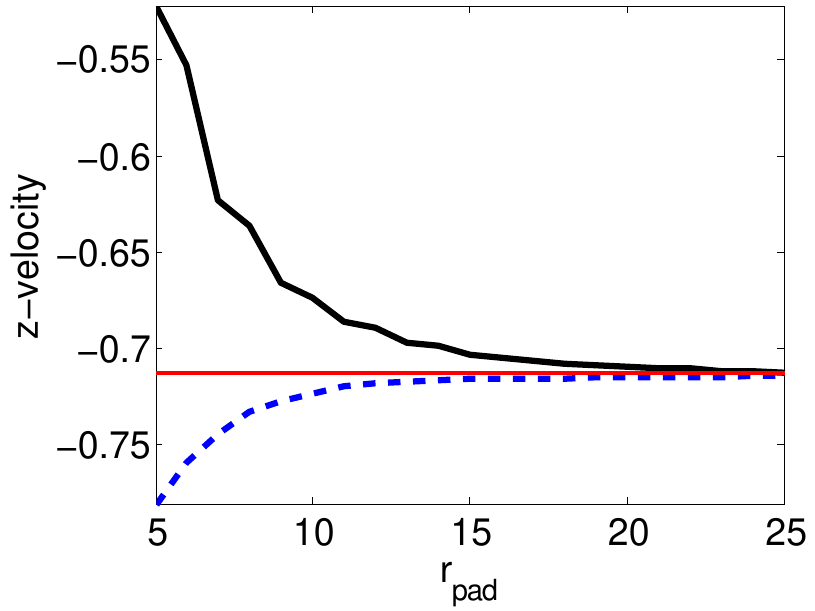}
\includegraphics[width=1.65in]{ 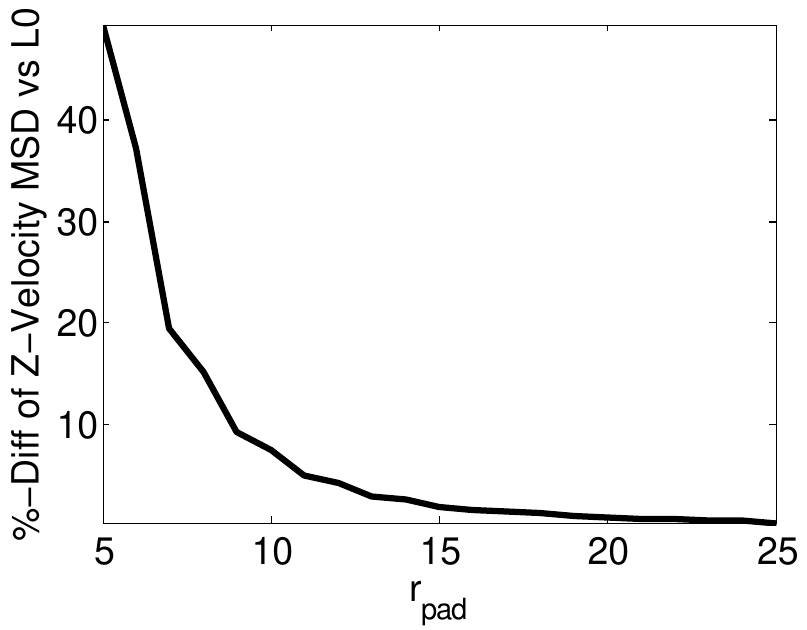}
\includegraphics[width=1.65in]{ 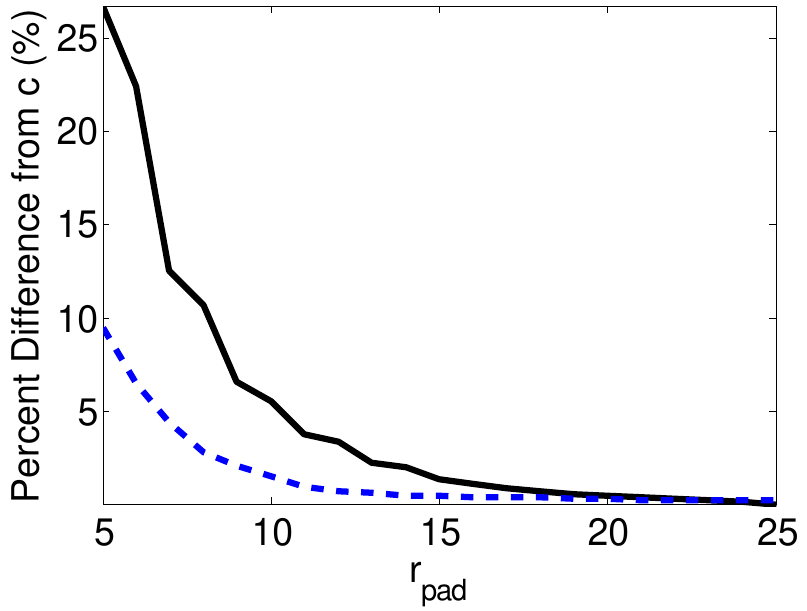}
\\
\includegraphics[width=1.65in]{ 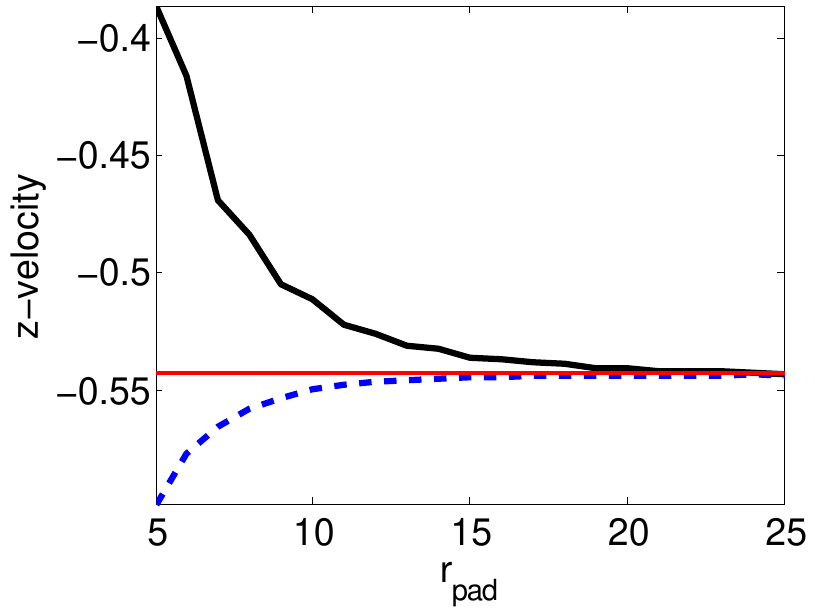}
\includegraphics[width=1.65in]{ 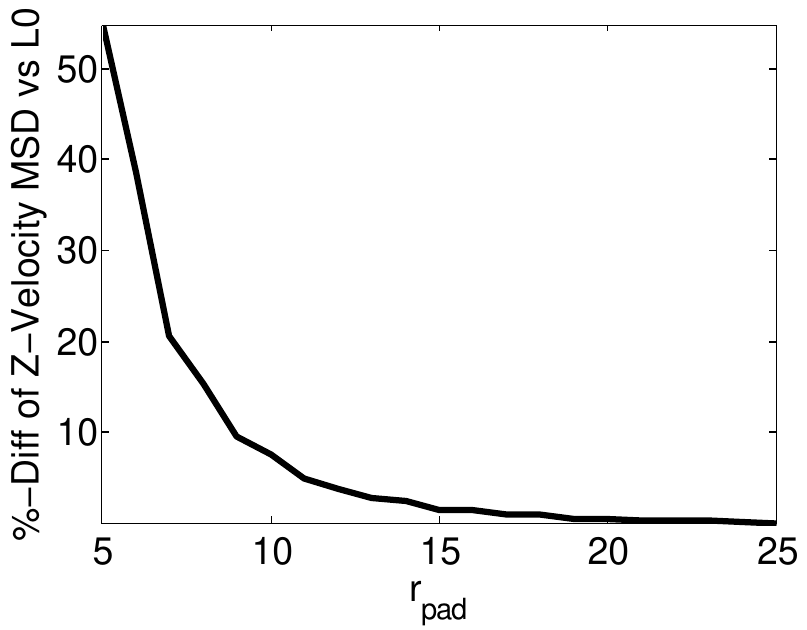}
\includegraphics[width=1.65in]{ 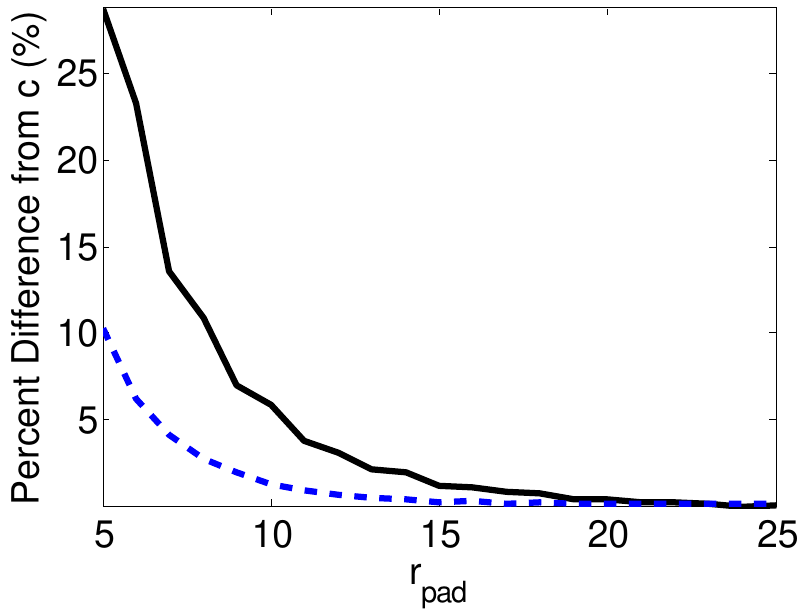}
\\
\includegraphics[width=1.65in]{ 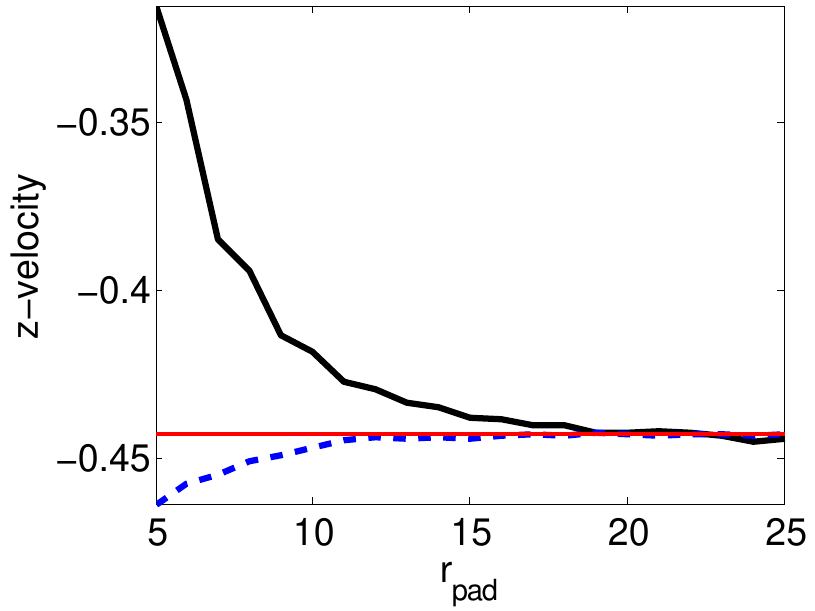}
\includegraphics[width=1.65in]{ 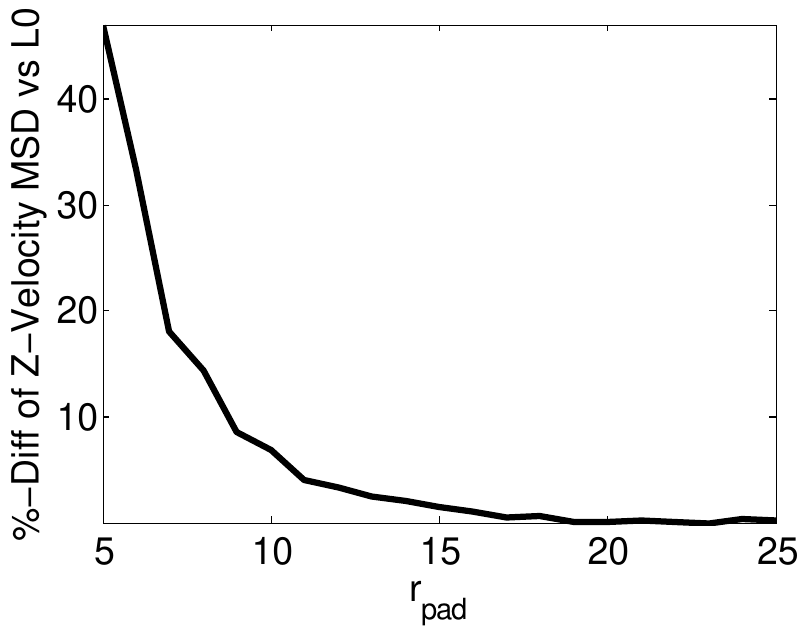}
\includegraphics[width=1.65in]{ 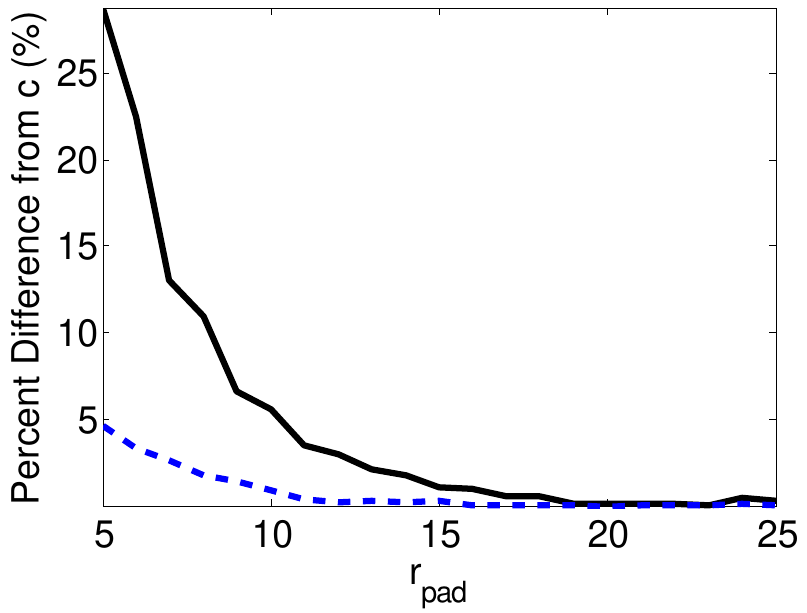}
\\
\includegraphics[height=2in]{ 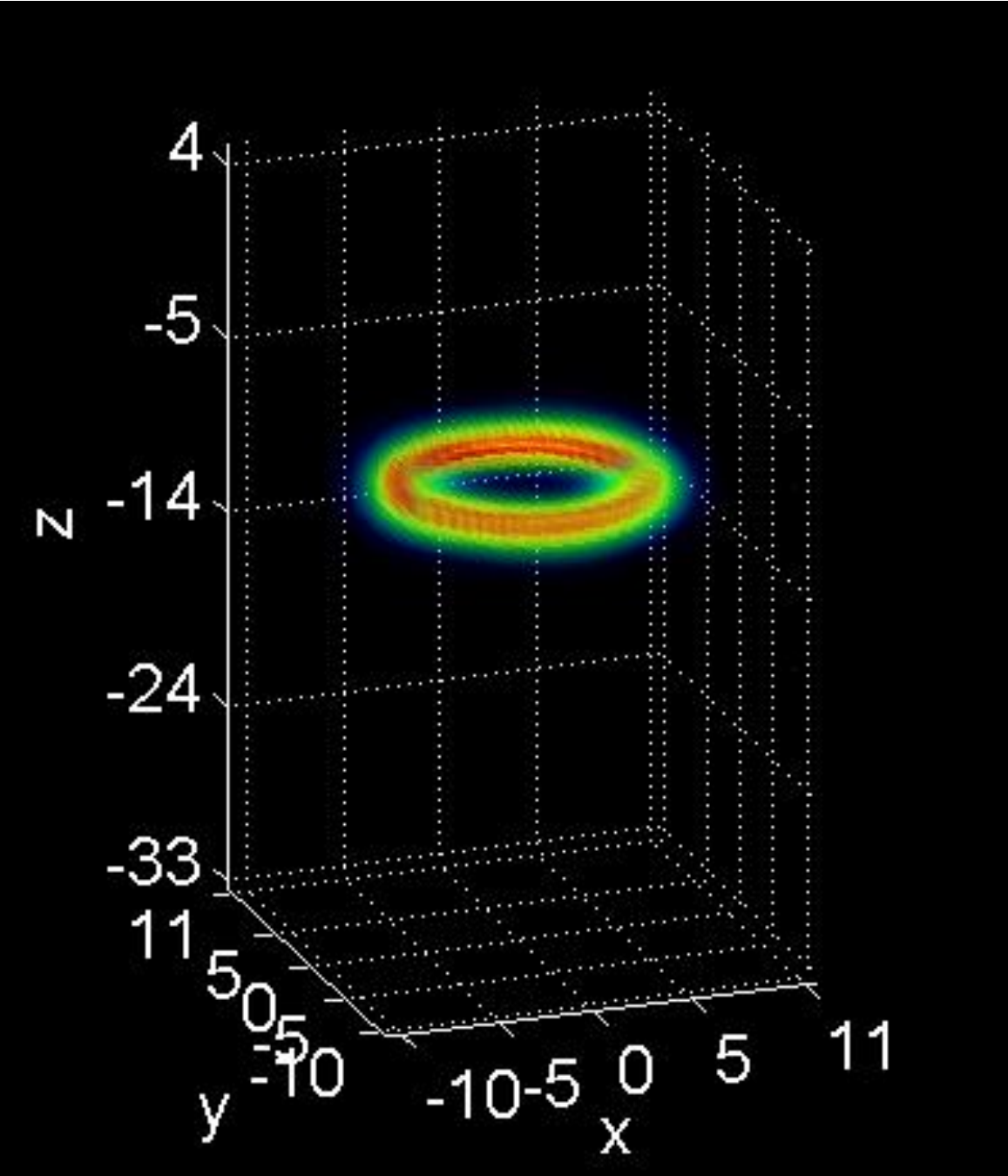}
\includegraphics[height=2in]{ 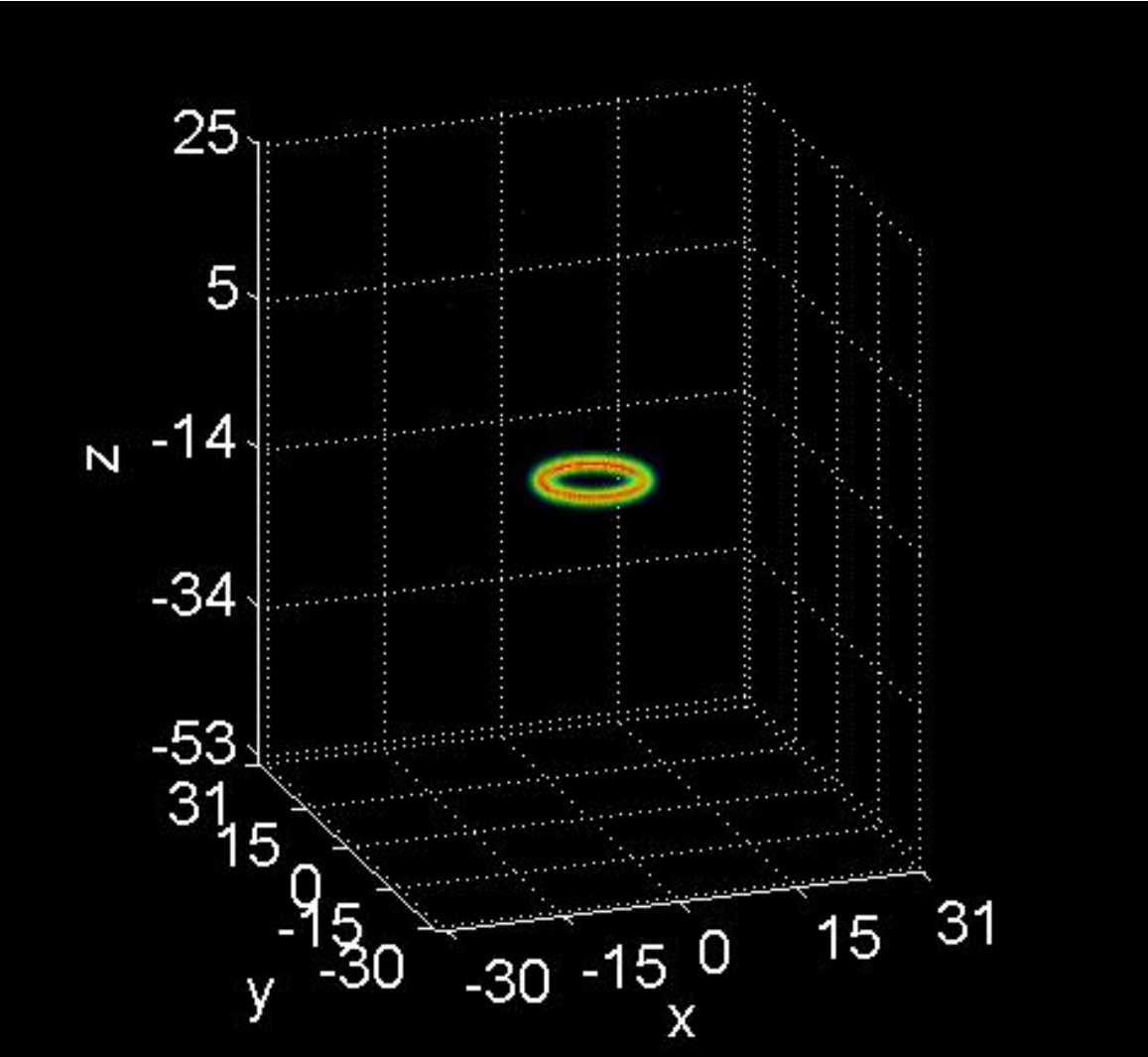}
\caption[MSD boundary condition vortex ring velocity tests.]{(Color online) Comparison of the velocity of a dark vortex ring in the NLSE over various grid sizes using the MSD and L0 boundary conditions.  Top-left row: The observed velocity of the vortex ring versus $r_{\mbox{\scriptsize pad}}\in[5,25]$ (the minimum allowed distance distance from the vortex core to the grid boundary allowing for the travel distance) for the MSD (thick solid line), L0 (blue dashed line) boundary conditions (the predicted value of the velocity (thin red line) is shown for comparison. Top-middle row: The percent-difference between the velocities observed using the MSD and L0 boundary conditions. Top-right row: The percent differences of the velocities observed compared to the predicted velocity.  Top-to-bottom: Results for a vortex ring of radius $4$, $6$, and $8$.   Bottom row: Volumetric rendering of the vortex ring at $t=20$ using the MSD boundary condition for $r_{\mbox{\scriptsize pad}}=5$ (left) and $r_{\mbox{\scriptsize pad}}=25$ (right).  The vortex rings are simulated to an end-time of $t=50$ using a time-step size of $k=0.04$ and a spatial-step size of $h=1/2$.\label{f:msdvrradtests}}
\end{figure}

In Fig.~\ref{f:msdvrradtests}, we show the results of simulating vortex rings of radii $4$, $6$, and $8$ for $50$ time units for grid sizes determined by a maximum distance (denoted $r_{\mbox{\scriptsize pad}}$) from the core of the vortex ring in each direction (taking the travel distance in the $z$-direction into account).  The difference between the velocity of the vortex rings when using the MSD boundary condition is compared to that of using the L0, as well as how each compare to the predicted value recomputed using the method of Ref.~\cite{VR_NLSE_VEL_71} (for details on the generation of the vortex rings and the tracking procedure, see Ref.~\cite{ME_DISS}).
We see that both the MSD and L0 boundary conditions converge to a constant velocity very close to that of the predicted value.  It is also noteworthy that the convergence happens at roughly the same distance from the vortex ring independent of the vortex ring's radius.

\begin{figure}[htb!]
\centering
\includegraphics[width=0.95in]{ 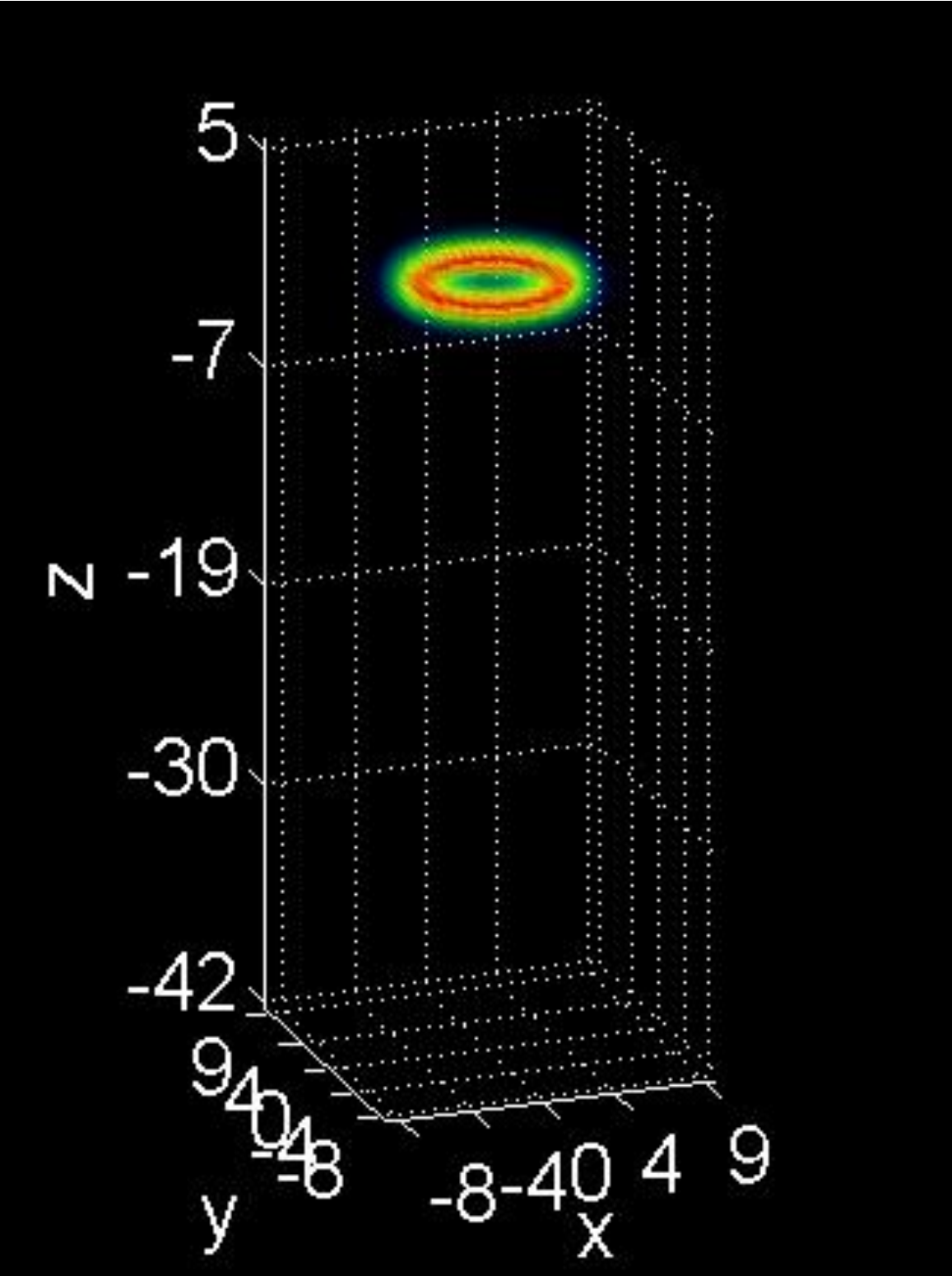}
\includegraphics[width=0.95in]{ 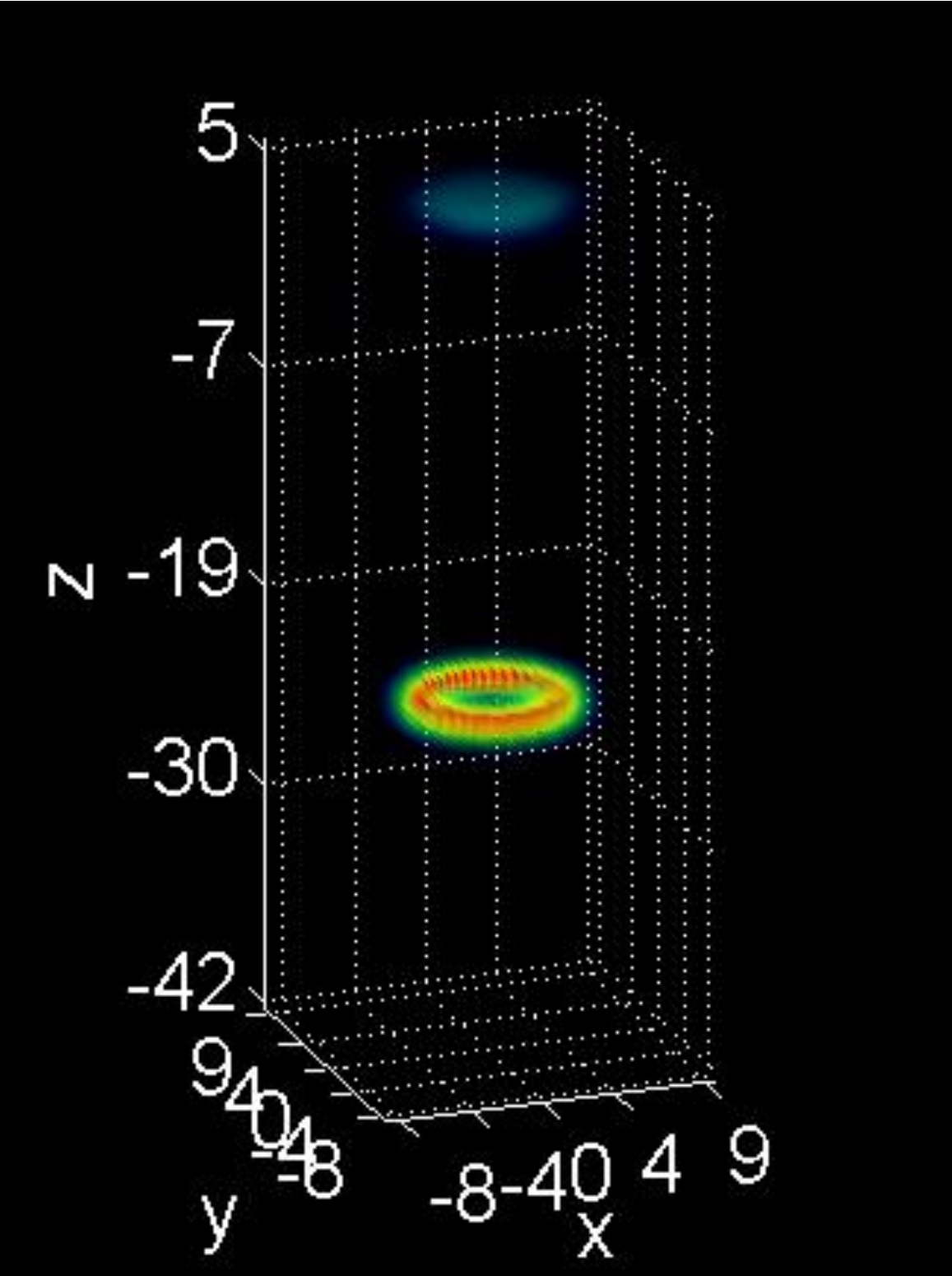}
\includegraphics[width=0.95in]{ 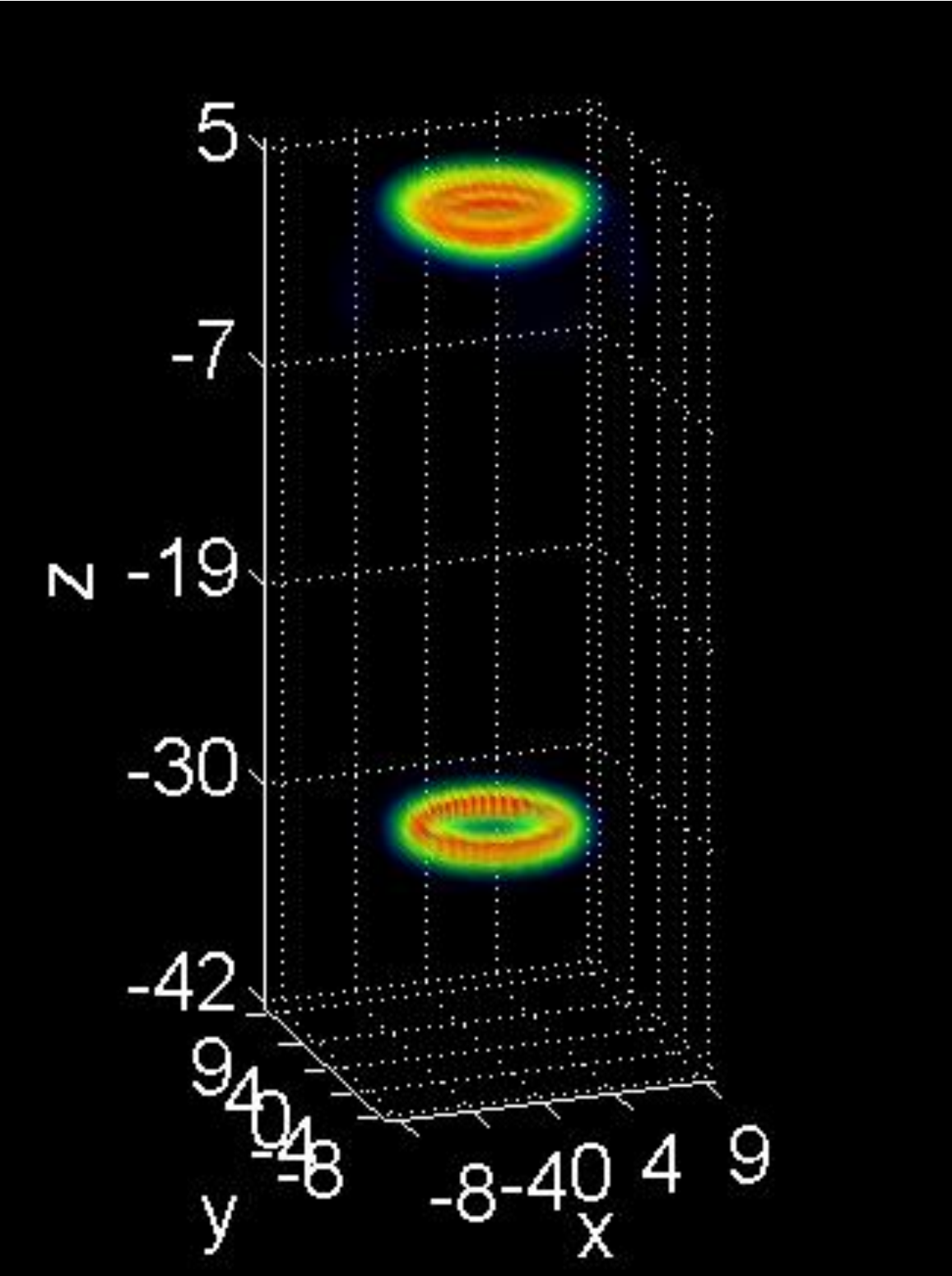}
\includegraphics[width=0.95in]{ 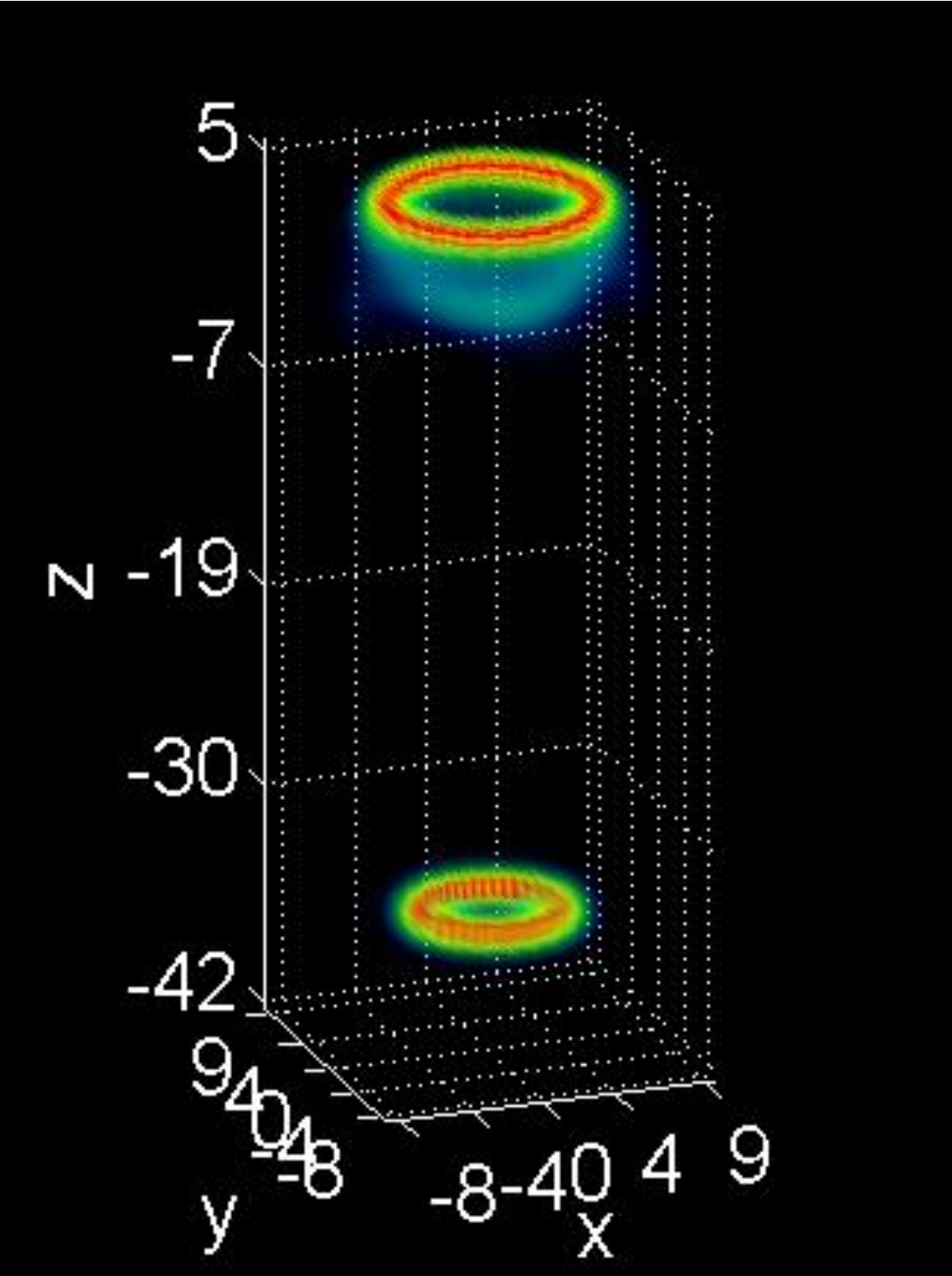}
\includegraphics[width=0.95in]{ 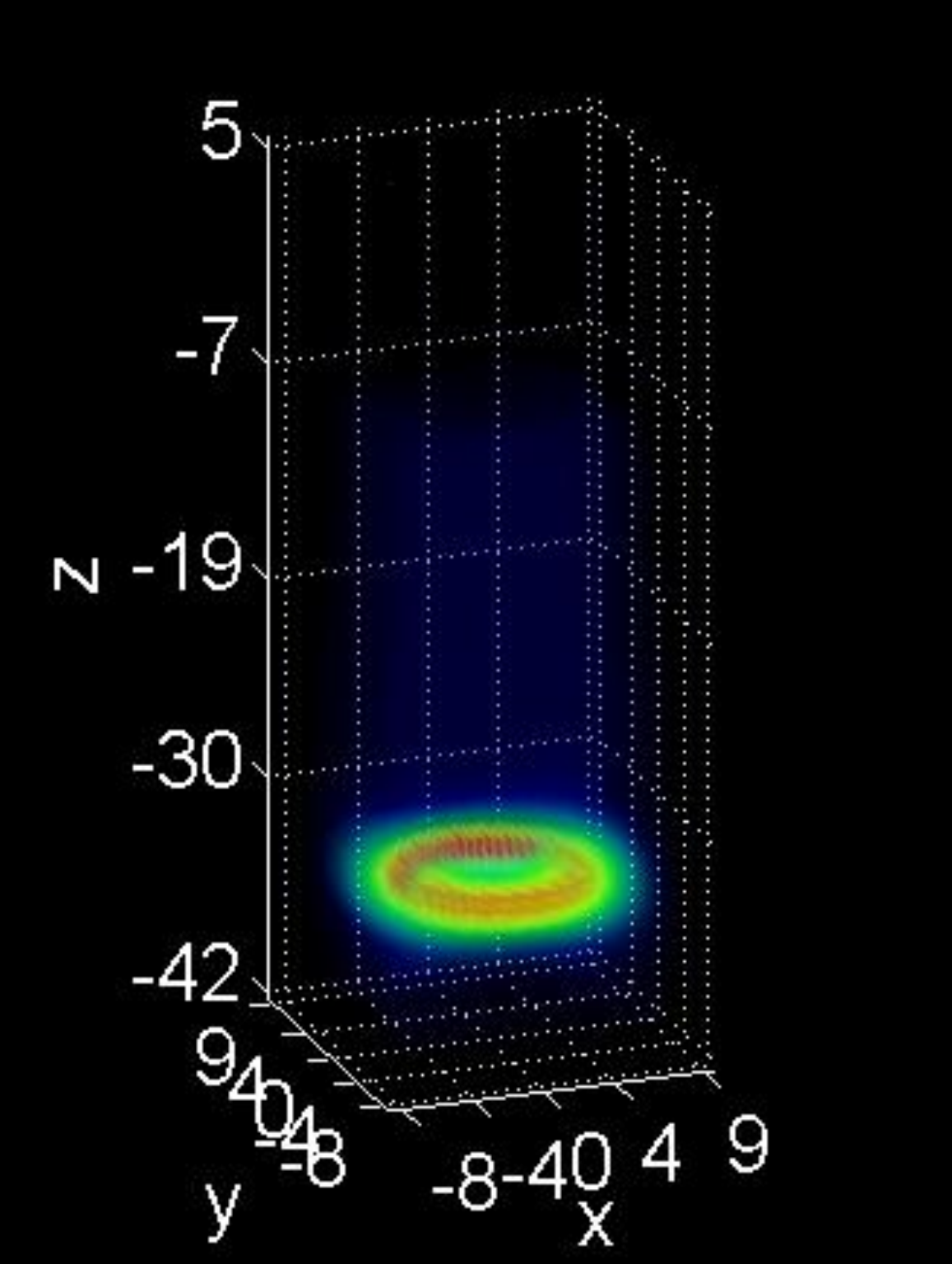}
\caption[Spontaneous vortex ring generation using L0 BC.]{(Color online) Example of spurious vortex ring generation when using L0 boundary condition on a small grid.  The snapshots of the vortex ring are shown for $r_{\mbox{\scriptsize pad}}=5$ at times $t=0$, $t=30$, $t=39$, and $t=45$.  The far right image is the same vortex ring at time $t=65$ using the MSD boundary condition (a later time was needed to simulate the vortex ring to the same grid position due to the slower velocity of the ring when using the MSD boundary condition as seen in Fig.~\ref{f:msdvrradtests}).  The numerical parameters used are the same as described in Fig.~\ref{f:msdvrradtests}.\label{f:vrl0err}}
\end{figure}

From Fig.~\ref{f:msdvrradtests} we notice that the L0 boundary condition converges to a constant velocity more rapidly than the MSD.  However, the L0 boundary condition can suffer from additional difficulties with smaller grid sizes, namely the spontaneous creation of spurious vortex rings due to the L0 boundary condition's problem in maintaining a proper phase structure as was shown in Fig.~\ref{f:2d1vort}.  An example of this is shown in Fig~\ref{f:vrl0err}.
We see that at such a low $r_{\mbox{\scriptsize pad}}$ value, the MSD boundary condition causes the vortex ring to become more distorted when compared to the L0 boundary condition (however, this should be weighed considering the fact that the vortex ring was simulated for a full $20$ time units longer than the case of the L0 boundary condition).  More importantly, the L0 boundary condition creates a spurious vortex ring near the top boundary of the grid.  Although one would not use such a low $r_{\mbox{\scriptsize pad}}$ value in either case, the error induced in the phase by the L0 boundary condition may greatly effect the dynamics of the vortex rings, especially for long simulations.

Lastly, let us comment upon the stability of the numerical methods that we used.
It is well known that boundary conditions can adversely effect the stability of
numerical simulations to the point where an otherwise stable scheme can become
unstable \cite{FD_PDE_BOOK}.  Therefore, it is necessary to consider the
stability effects of the MSD boundary condition.
However, since the MSD boundary condition is general in nature, its effect on
the stability of the simulations will depend on the governing PDE being used,
as well as the form of the overall numerical scheme being implemented.
Therefore, no general statements about the MSD boundary condition's stability
effects can be made.
That being said, in the specific case of the NLSE using the RK4 time-stepping
with the CD or 2SHOC spatial schemes that we have used for the examples in
Sec.~\ref{s:num}, the stability effects of the MSD boundary condition can be
addressed.  An analysis of the stability effects of using the MSD
(and the L0) boundary conditions has been carried out as part of our prior
study on the stability of RK4 schemes applied to the NLSE in Ref.~\cite{ME_RK4STB} where it was found that the MSD boundary condition did not have a significant effect on the stability bounds of the overall scheme.  Although those stability results are in no way general, as they are specific to the
NLSE using the RK4+CD and RK4+2SHOC schemes, they may be indicative
of similar results with difference governing PDEs and methods.

\section{Conclusions}
\label{s:con}
We have shown the formulation of a modulus-squared Dirichlet (MSD) boundary condition for numerical simulations of time-dependent complex partial differential equations. The standard form of the MSD boundary condition is given as a boundary value of the time-derivative of the solution as a function of the solution at that point, as well as the solution and its time-derivative at the closest interior point.  It is easily expressed as
\[
\left.\frac{\partial \Psi}{\partial t}\right|_b \approx i\,\mbox{Im}\left[  \frac{1}{\Psi_{b-1}}  \left.\frac{\partial \Psi}{\partial t}\right|_{b-1}  \right]\,\Psi_b,
\]
where the subscripts $b$ and $b-1$ refer to a boundary point and the closest interior point to the boundary respectively.

Through multidimensional numerical examples of the MSD boundary condition applied to the nonlinear {Schr{\"o}dinger} equation, we have shown that it is extremely effective in terms of non-interference with internal dynamics, as well as in requiring smaller grid sizes when compared to other boundary conditions currently in use (such as setting the Laplacian to zero at the boundaries).  This is especially true in simulations of coherent structures which exhibit vorticity (such as vortices and vortex rings) and in simulations of solutions which have a co-moving background velocity associated with them.  In the latter case, the MSD boundary condition seems to be the {\em only} simple boundary condition which can handle such co-moving back-flows.  We conclude that the MSD boundary condition applied to complex-valued problems is as effective as the standard Dirichlet boundary condition for problems with a constant value at the boundary, and nearly as easy to implement.

%
%

\section*{Acknowledgments}
We are thankful to an anonymous referee for providing insight into 
the continuous formulation of the MSD boundary condition which was
very helpful in clarifying the overall derivation given in Sec.~\ref{s:derive}.
This research was supported by NSF-DMS-0806762 and the Computational
Science Research Center at San Diego State University.

\def\myitemsep{5pt}
\bibliographystyle{mio}
\bibliography{MSD}  
\end{document}